\documentclass[oneside,english]{amsart}
\usepackage[LGR,T1]{fontenc}
\usepackage[latin9]{inputenc}
\usepackage{geometry}
\usepackage{units}
\usepackage{mathrsfs}
\usepackage{mathtools}
\usepackage{amsbsy}
\usepackage{amstext}
\usepackage{amsthm}
\usepackage{amssymb}
\usepackage{stmaryrd}
\usepackage{stackrel}
\usepackage{graphicx}

\makeatletter

\newcommand{\noun}[1]{\textsc{#1}}
\DeclareRobustCommand{\greektext}{%
  \fontencoding{LGR}\selectfont\def\encodingdefault{LGR}}
\DeclareRobustCommand{\textgreek}[1]{\leavevmode{\greektext #1}}
\ProvideTextCommand{\~}{LGR}[1]{\char126#1}

\providecommand{\tabularnewline}{\\}

\numberwithin{equation}{section}
\numberwithin{figure}{section}
\theoremstyle{plain}
\newtheorem{thm}{\protect\theoremname}
\theoremstyle{remark}
\newtheorem{rem}[thm]{\protect\remarkname}
\theoremstyle{plain}
\newtheorem{prop}[thm]{\protect\propositionname}

\usepackage{babel}
\usepackage{bbold}
\usepackage[all]{xy}

\makeatother

\usepackage{babel}
\providecommand{\propositionname}{Proposition}
\providecommand{\remarkname}{Remark}
\providecommand{\theoremname}{Theorem}

\begin{document}
\title{collineation groups of octonionic and split-octonionic planes}
\author{Daniele Corradetti$^{*}$, Alessio Marrani$^{\dagger}$, Francesco
Zucconi$^{\ddagger}$}
\begin{abstract}
We present a Veronese formulation of the octonionic and split-octonionic
projective and hyperbolic planes. This formulation of the incidence
planes highlights the relationship between the Veronese vectors and
the rank-1 elements of the Albert algebras over octonions and split-octonions,
yielding to a clear formulation of the relationship with the real
forms of the Lie groups arising as collineation groups of these planes.
The Veronesean representation also provides a novel and minimal construction
of the same octonionic and split-octonionic planes, by exploiting
two symmetric composition algebras: the Okubo algebra and the paraoctonionic
algebra. Besides the intrinsic mathematical relevance of this construction
of the real forms of the Cayley-Moufang plane, we expect this approach
to have implications in all mathematical physics related with exceptional
Lie Groups of type $G_{2},F_{4}$ and $E_{6}$.
\end{abstract}

\maketitle
\tableofcontents{}

\section{Introduction\label{sec:Introduction} }

Recent advances in the algebraic formulation of the Standard Model
of elementary particle physics \cite{Ba21,FH22,MDW,Ma21,To23} attracted
new attention within the mathematical physics community which study
octonions. However, while the algebra of octonions garnered some attention,
the algebra of split-octonions, which is not a division algebra, largely
remained on the periphery. Nevertheless, a recent work of Penrose
\cite{Pen22} highlighted that quantized bi-twistors possess a hitherto
unobserved $G_{2\left(2\right)}$ structure, which directly associates
them with split-octonions. This hints to the relevance of the real
forms of Lie algebras which relate to split octonions themselves \cite{CMR}.
Among these, the real forms of Lie algebras of type $E_{6}$ (and
the corresponding Lie groups) \cite{CCMA-Magic}, have a remarkable
relevance, also in relation to Grand Unified Theories \cite{BH10}.
Since groups of type $F_{4}$ and $E_{6}$ can respectively be realized
as automorphism groups and reduced structure groups of Albert algebras
over octonions or split octonions, it is clear that such algebras
have attracted new attention.

A geometric construction of the automorphism (resp. reduced structure
group) of the Albert algebra $\mathfrak{J}_{3}\left(\mathbb{O}\right)$
is achieved through the isometry (resp. collineation)group of the
corresponding octonionic projective plane (aka Cayley-Moufang plane)
$\mathbb{O}P^{2}$. This motivates us to present an useful and constructive
way to introduce the octonionic and split-octonionic projective and
hyperbolic planes, i.e. $\mathbb{O}P^{2}$, $\mathbb{O}H^{2}$, $\mathbb{O}_{s}P^{2}$
and $\mathbb{O}_{s}H^{2}$, through the use of a special set of vectors
called ``Veronese vectors'' \cite[Sec. 1.16]{Compact Projective}.
It should however be remarked that the latter two incidence planes,
introduced by Springer and Veldkamp in \cite{SpVa68}, are projective
and hyperbolic planes only in a vague sense since they do not satisfy
all projective axioms because split-octonions $\mathbb{O}_{s}$ are
not a division algebra. Namely, the number of points incident with
two unique lines might exceed one, akin to the number of lines passing
through two distinct points.\medskip{}

The plan of the paper is as follows. After defining the octonionic
planes in Sec. \ref{sec:Veronese-formulation-of}, we focus on their
relationships with the corresponding real forms of the Albert algebra,
i.e. $\mathfrak{J}_{3}\left(\mathbb{O}\right)$, $\mathfrak{J}_{2,1}\left(\mathbb{O}\right)$,
$\mathfrak{J}_{3}\left(\mathbb{O}_{s}\right)$ and $\mathfrak{J}_{2,1}\left(\mathbb{O}\right)$,
whose Veronese vectors are the rank-1 elements. We then study groups
of motions over the octonionic projective plane recovering $\text{E}_{6\left(-26\right)}$
and $\text{F}_{4\left(-52\right)}$ as the collineation resp. isometry
group of $\mathbb{O}P^{2}$. With a similar treatment, we recover
other real forms of $\text{E}_{6}$ and $\text{F}_{4}$: namely, $\text{E}_{6\left(6\right)}$
and $\text{F}_{4\left(4\right)}$ for the split case $\mathbb{O}_{s}P^{2}$;
$\text{E}_{6\left(-26\right)}$ and $\text{F}_{4\left(-20\right)}$
for the octonionic hyperbolic plane $\mathbb{O}H^{2}$; and again
$\text{E}_{6\left(6\right)}$ and $\text{F}_{4\left(4\right)}$ in
the split-hyperbolic case $\mathbb{O}_{s}H^{2}$. If one focus on
the collineation group of the projective planes that fixes non-degenerate
quadrangles, then one recovers $\text{G}_{2\left(-14\right)}$ and
$\text{G}_{2\left(2\right)}$ for the octonionic and split-octonionic
cases respectively \cite{Compact Projective,corr Notes Octo,Corr RealF}.

In Sec. \ref{sec:Realisations-through-symmetric} we then present
a novel, in a sense minimal, construction for obtaining the same planes
making use of symmetric composition algebras instead of Hurwitz algebras.
This is a natural continuation of our previous work on Okubo algebras.
Indeed, with a slight modification of Veronese conditions, the paraoctonionic
and the Okubo algebras are suitable for the construction of the Cayley
plane \cite{CMZ}. In this work we extend the result, showing that
the same mathematical setup can be used for the construction of all
previous incidence planes using symmetric composition algebras instead
of Hurwitz. Finally, in Sec. \ref{sec:Complex-Cayley-plane} we show
how this setup can be used for also on the field of complex numbers
$\mathbb{C}$, thus obtaining the complex groups $\text{E}_{6}^{\mathbb{C}},\text{F}_{4}^{\mathbb{C}}$
and $\text{G}_{2}^{\mathbb{C}}$.

\section{\label{sec:octonions-and-split-octonions}Octonions and split-octonions }

An\emph{ algebra} is a vector space $A$ over a field $\mathbb{F}$
(that we will assume to be $\mathbb{R}$ and $\mathbb{C}$) with a
bilinear multiplication. The algebra $A$ is said to be\emph{ commutative}
if $x\cdot y=y\cdot x$ for every $x,y\in A$; \emph{associative}
if satisfies $x\cdot\left(y\cdot z\right)=\left(x\cdot y\right)\cdot z$;
\emph{alternative} if $x\cdot\left(y\cdot y\right)=\left(x\cdot y\right)\cdot y$;
and finally, \emph{flexible} if $x\cdot\left(y\cdot x\right)=\left(x\cdot y\right)\cdot x$.
Furthermore, if the algebra has a non degenerate norm $n$ that upholds
the multiplicative property

\begin{align}
n\left(x\cdot y\right) & =n\left(x\right)n\left(y\right),\label{eq:comp(Def)}
\end{align}
 for every $x,y\in A$, is called a \emph{composition} algebra \cite{ElDuque Comp}
and is denoted with the triple $\left(A,\cdot,n\right)$ or simply
as $A$ if there are no reason for ambiguity.

To qualify as an algebra, $A$ must be a group over addition $+$.
However, there are no similar group requirements for the bilinear
product, allowing it to be non-associative or lack an identity element.
Algebras with an element $1$ such that $1\cdot x=x\cdot1=x$ are
called \emph{unital}.

According to the Hurwitz theorem, only two eight-dimensional composition
algebras that are also unital exist: the octonions $\mathbb{O}$ and
split-octonions $\mathbb{O}_{s}$. Both of them are non-associative.
In fact, they almost share the same algebraic property with the exception
of a different quadratic form signature. As consequence, the former,
i.e. the octonions $\mathbb{O}$ are a division algebra, while the
latter, i.e. the split-octonions $\mathbb{O}_{s}$ are not. In Sec.
\ref{sec:Realisations-through-symmetric} we will see that dropping
the requirement of the algebra being unital other three 8-dimensional
algebras appear: the paraoctonions $p\mathbb{O}$, the split-paraoctonions
$p\mathbb{O}_{s}$ and the real Okubo algebra $\mathcal{O}$.

One can directly and abstractly define both octonions and split-octonions
by considering the complex vector space $\mathbb{O}_{\mathbb{C}}$
endowed with the multiplication table in Tab. \ref{tab:Split Octonions-1-1}
over the base $\left\{ e_{1},e_{2},u_{1},u_{2},u_{3},v_{1},v_{2},v_{3}\right\} $
obtained starting from two idempotents $e_{1}$ and $e_{2}$ and two
three-dimensional real vector subspace such that $V=\left\{ e_{1}\cdot v=0\right\} $
and $U=\left\{ e_{1}\cdot u=0\right\} $. The resulting algebra is
that of the complex octonions \cite{ZSSS}.
\begin{figure}
\centering{}\includegraphics[scale=0.07]{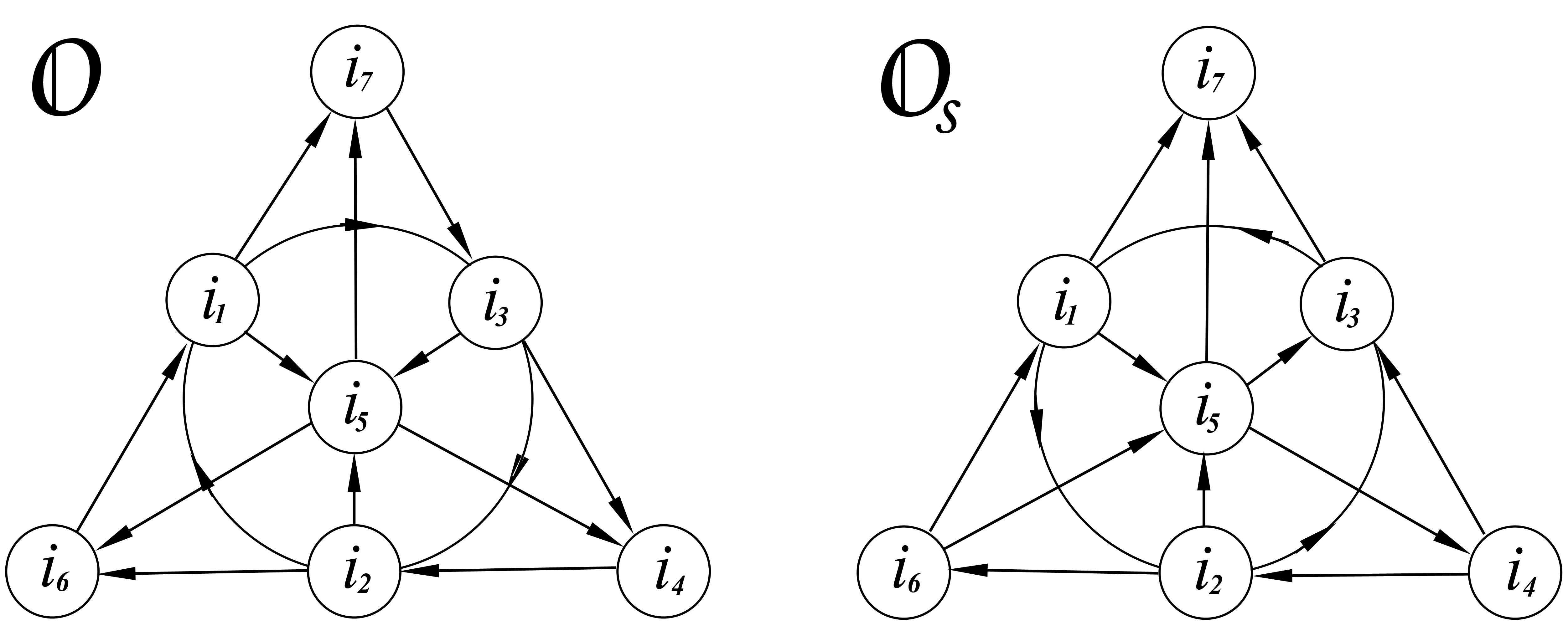}\caption{\label{fig:Octo-SplitOcto}Multiplication rule of octonions $\mathbb{O}$
and of split-octonions $\mathbb{O}_{s}$ as real vector space defined
over the basis $\left\{ \text{i}_{0}=1,\text{i}_{1},\text{i}_{2},\text{i}_{3},\text{i}_{4},\text{i}_{5},\text{i}_{6},\text{i}_{7}\right\} $.
Lines in the Fano plane identify associative triples of the product
and the arrow indicates the sign (positive in the sense of the arrow
and negative in the opposite sense). In addition to the previous rules
it is intended that $\text{i}_{k}^{2}=-1$ in the case of octonions
$\mathbb{O}$ and $\text{i}_{k}^{2}=-1$ for $k=1,2,3$ and $\text{i}_{k}^{2}=1$
otherwise.}
\end{figure}
 
\begin{table}
\begin{centering}
\begin{tabular}{|c|cc|ccc|ccc|}
\hline 
 & $e_{1}$ & $e_{2}$ & $u_{1}$ & $u_{2}$ & $u_{3}$ & $v_{1}$ & $v_{2}$ & $v_{3}$\tabularnewline
\hline 
$e_{1}$ & $e_{1}$ & 0 & $u_{1}$ & $u_{2}$ & $u_{3}$ & 0 & 0 & 0\tabularnewline
$e_{2}$ & 0 & $e_{2}$ & 0 & 0 & 0 & $v_{1}$ & $v_{2}$ & $v_{3}$\tabularnewline
\hline 
$u_{1}$ & 0 & $u_{1}$ & 0 & $v_{3}$ & $-v_{2}$ & $-e_{1}$ & 0 & 0\tabularnewline
$u_{2}$ & 0 & $u_{2}$ & $-v_{3}$ & 0 & $v_{1}$ & 0 & $-e_{1}$ & 0\tabularnewline
$u_{3}$ & 0 & $u_{3}$ & $v_{2}$ & $-v_{1}$ & 0 & 0 & 0 & $-e_{1}$\tabularnewline
\hline 
$v_{1}$ & $v_{1}$ & 0 & $-e_{2}$ & 0 & 0 & 0 & $u_{3}$ & $-u_{2}$\tabularnewline
$v_{2}$ & $v_{2}$ & 0 & 0 & $-e_{2}$ & 0 & $-u_{3}$ & 0 & $u_{1}$\tabularnewline
$v_{3}$ & $v_{3}$ & 0 & 0 & 0 & $-e_{2}$ & $u_{2}$ & $-u_{1}$ & 0\tabularnewline
\hline 
\end{tabular}
\par\end{centering}
\medskip{}

\caption{\label{tab:Split Octonions-1-1}Multiplication table of the complex
octonions in the canonical base from \cite{ZSSS}.}
\end{table}

Given the definition of complex octonions $\mathbb{O}_{\mathbb{C}}$,
as described, only two real algebras, when complexified, are isomorphic
to $\mathbb{O}_{\mathbb{C}}$: namely, the octonions $\mathbb{O}$
and the split-octonions $\mathbb{O}_{s}$. In other word these two
algebras are the only algebras such that 
\begin{equation}
\mathbb{O}_{\mathbb{C}}\cong\mathbb{C}\otimes\mathbb{O}\cong\mathbb{C}\otimes\mathbb{O}_{s}.
\end{equation}

A more practical way, to define both algebras is to consider the real
vector space over the basis $\left\{ \text{i}_{0}=1,\text{i}_{1},\text{i}_{2},\text{i}_{3},\text{i}_{4},\text{i}_{5},\text{i}_{6},\text{i}_{7}\right\} $
and define an appropriate multiplication for the octonions and the
split-octonions following that indicated by the Fano plane in Fig.
\ref{fig:Octo-SplitOcto}, with the addition of $\text{i}_{k}^{2}=-1$
in the case of octonions $\mathbb{O}$ and $\text{i}_{k}^{2}=-1$
for $k=1,2,3$ and $\text{i}_{k}^{2}=1$ for $k\neq1,2,3$ in the
case of split-octonions $\mathbb{O}_{s}$. As a consequence, given
the decomposition of an element $x\in\left\{ \mathbb{O},\mathbb{O}_{s}\right\} $
and its conjugate $\overline{x}$, i.e.,
\begin{align}
x & =x_{0}+\stackrel[k=1]{7}{\sum}x_{k}\text{i}_{k},\\
\overline{x} & =x_{0}-\stackrel[k=1]{7}{\sum}x_{k}i_{k},
\end{align}
the norm $n$ is assumes two different forms
\begin{align}
n\left(x\right) & =\left(x_{0}^{2}+x_{1}^{2}+x_{2}^{2}+x_{3}^{2}\right)+\left(x_{4}^{2}+x_{5}^{2}+x_{6}^{2}+x_{7}^{2}\right),\label{eq:octonionic norm-1}\\
n_{s}\left(x\right) & =\left(x_{0}^{2}+x_{1}^{2}+x_{2}^{2}+x_{3}^{2}\right)-\left(x_{4}^{2}+x_{5}^{2}+x_{6}^{2}+x_{7}^{2}\right).
\end{align}
Nevertheless, in both cases we have 
\begin{equation}
x\cdot\overline{x}=n\left(x\right)1,
\end{equation}
as always is in unital composition algebras. 

The most interesting feature of split-octonions is their suitability
in giving rise to topological variants of the same geometrical objects
arising from octonionic constructions. For our purposes, it is very
important to notice that the automorphism group of $\mathbb{O}_{s}$
is still an exceptional Lie group of type $G_{2}$, but unlike the
octonionic case, is not the compact one $G_{2\left(-14\right)}$ but
the non-compact one $G_{2\left(2\right)}$. 

\section{\label{sec:Albert-algebras}Albert algebras}

Jordan algebras are commutative and non-associative algebras that
were originally introduced to model the foundations of quantum mechanics
\cite{Jordan}. Although they did not achieve this intended purpose,
they turned out to have numerous and deep connections with many branches
of mathematics, such as Lie theory, symmetric spaces, and projective
geometry. In fact, the investigation of Jordan algebras---especially
in their relationships with exceptional Lie groups and associated
algebras---represents more than three decades of impactful mathematical
research from the 20th century. A comprehensive summary goes beyond
the scope of just a few pages, but readers can consult \cite{McCr}
for historical insights. 

The goal of this section is to is to provide readers with the essential
theoretical background necessary for subsequent sections. We will
review the definition of Jordan algebras as proposed by Springer \cite{Krut,McCr},
which originates from a cubic norm on a vector space. Our emphasis
on this particular definition is motivated by its explicit connection,
via the map $\#$, to Veronese conditions, a pivotal topic in the
ensuing sections. We will then delve into the study of exceptional
Jordan algebras, commonly known as \emph{Albert algebras}, i.e., $\mathfrak{J}_{3}\left(\mathbb{O}\right)$,
$\mathfrak{J}_{3}\left(\mathbb{O}_{s}\right)$, $\mathfrak{J}_{2,1}\left(\mathbb{O}\right)$
and $\mathfrak{J}_{2,1}\left(\mathbb{O}_{s}\right)$, that we analyse
within a unified framework. Finlly, we analyse the automorphisms of
the Albert algebras and cubic-norm-preserving transformations (also
named reduced structure group), showing how exceptional Lie groups
of type $\text{F}_{4}$ and $\text{E}_{6}$ naturally emerge in this
context, as originally demonstrated by Chevalley and Schafer in \cite{ChSch}.
This discovery warranted further inquiries into the ties between Jordan
algebras and Lie groups. In fact, at the present the easiest realisation
of the real forms of such groups is by means of Albert algebras \cite{Vinberg}.
It is worth saying that also $E_{7}$ and even $E_{8}$ can be recovered
within the same framework, as shown by Freudenthal in \cite{Fr54}.
An important reference on Jordan algebras is McCrimmon's book \cite{McCr},
along with Jacobson's treatise \cite{Jac}; additional key references
include Schafer \cite{Scha} and the collaborative work of Springer
and Valdekamp \cite{Veldkamp}. For the relation between Jordan algebras
and Lie groups, we will often refer to Yokota \cite{Yokota} which
is one of the most valuable resources. 

\subsection{\label{sec:Jordan-algebras}Jordan algebras}

A Jordan algebra $\mathfrak{J}$ is a vector space endowed with a
commutative bilinear product that satisfies the Jordan identity

\begin{equation}
\left(x\circ x\right)\circ\left(x\circ y\right)=x\circ\left(x\circ\left(x\circ y\right)\right),
\end{equation}
for every $x,y\in\mathfrak{J}$. By definition, Jordan algebras are
commutative and thus flexible; however, they generally they are not
associative. The prototypical example of a Jordan algebra is derived
from an associative algebra $\mathfrak{A}$ defining the Jordanian
product as
\begin{equation}
x\circ y=\frac{1}{2}\left(xy+yx\right).
\end{equation}
Jordan algebras that do not originate from an underlying associative
algebra are termed \emph{exceptional}.

A large class of Jordan algebras, arise from a cubic form $N$ called
\emph{Jordan admissible} . These are called \emph{cubic Jordan algebras}
and are due to Springer \cite{Springer}. Comprehensive reviews on
this topic are available in \cite{McCr} and in \cite{Krut}. Let
$N$ be a polynomial map over a real vector space $V$ such that $N$
is cubic (homogeneous of degree three), i.e. 
\begin{equation}
N\left(\lambda x\right)=\lambda^{3}N\left(x\right)
\end{equation}
for every $\lambda\in\mathbb{R}$ and $x\in V$ and the full linearisation
of $N\left(x\right)$ is a trilinear map given given by
\begin{equation}
N\left(x,y,z\right)=\frac{1}{6}\left(N\left(x+y+z\right)-N\left(x+y\right)-N\left(x+z\right)-N\left(y+z\right)+N\left(x\right)+N\left(y\right)+N\left(z\right)\right),
\end{equation}
with $x,y,z\in V$. Moreover, let $c\in V$ be a basepoint for $N$,
i.e. $N\left(c\right)=1$, and define a linear map called \emph{trace}
\begin{equation}
\text{Tr}\left(x\right)=3N\left(c,c,x\right),
\end{equation}
a quadratic map, which is the usual \emph{quadratic norm}, 
\begin{equation}
S\left(x\right)=3N\left(x,x,c\right),
\end{equation}
a bilinear map, which is the \emph{polar form} of the previous quadratic
norm,
\begin{equation}
S\left(x,y\right)=6N\left(x,y,c\right),
\end{equation}
 and finally a trace-bilinear form 
\begin{equation}
\left(x,y\right)=\text{Tr}\left(x\right)\text{Tr}\left(y\right)-S\left(x,y\right)
\end{equation}
for every $x,y\in V$. Then, the cubic norm $N$ is called \emph{Jordan
admissible} if the trace bilinear form is non-degenerate at the base
point $c$ and the quadratic \emph{sharp} map defined\footnote{Notice that this very definition implies that $x^{\#}$ is homogeneous
of degree two in x, namely $\left(\lambda x\right)^{\#}=\lambda^{2}x^{\#}$.} uniquely as $\left(x^{\#},y\right)=3N\left(x,x,y\right)$ satisfies
the adjoint identity 
\begin{equation}
\left(x^{\#}\right)^{\#}=N\left(x\right)x.
\end{equation}
A vector space with a Jordan admissible cubic norm can be endowed
with a Jordan structure as the following theorem assures
\begin{thm}
\emph{\cite[Section I.3.8]{McCr}} Any vector space with an admissible
cubic form can be converted into a Jordan algebra with unit $c=\text{\textbf{1}}$
and Jordan product given by
\begin{equation}
x\circ y=\frac{1}{2}\left(x\times y+\text{Tr}\left(x\right)y+\text{Tr}\left(y\right)x-S\left(x,y\right)\text{\textbf{1}}\right),
\end{equation}
where $x\times y$ is twice the linearization of the sharp map, i.e.
$x\times y=\left(x+y\right)^{\#}-x^{\#}-y^{\#}$. Additionally, every
element of the Jordan algebra satisfy the cubic polynomial
\begin{equation}
x^{3}-\text{Tr}\left(x\right)x^{2}+S\left(x\right)x-N\left(x\right)\text{\textbf{1}}=0
\end{equation}
 and, also, $x^{\#}=x^{2}-\text{Tr}\left(x\right)x+S\left(x\right)\text{\textbf{1}}$.
\end{thm}

This foundational mathematical structure is conveniently applied to
vector spaces without invoking Hermitian matrices or the concept of
conjugation. Nevertheless, once we apply this construction to three
by three Hermitian matrices over unital composition algebras, i.e.
$H_{3}\left(\mathbb{K}\right)$ where $\mathbb{K}$ is an Hurwitz
algebra with norm $n$, we obtain the classical rank-three Jordan
algebras $\mathfrak{J}_{3}\left(\mathbb{K}\right)$. Indeed, let $X$
be an arbitrary element of $H_{3}\left(\mathbb{K}\right)$ with form
\begin{equation}
X=\left(\begin{array}{ccc}
\lambda_{1} & x_{3} & \overline{x}_{2}\\
\overline{x}_{3} & \lambda_{2} & x_{1}\\
x_{2} & \overline{x}_{1} & \lambda_{3}
\end{array}\right),\label{eq:Hermitian element}
\end{equation}
with $\lambda_{1},\lambda_{2},\lambda_{3}\in\mathbb{R}$, $x_{1},x_{2},x_{3}\in\mathbb{K}$
and $\overline{x}$ is the conjugate in $\mathbb{K}$ under the canonical
conjugation, i.e. $\overline{x}=\left\langle x,1\right\rangle 1-x$.
Then, let the basepoint $c$ be the identity matrix and consider the
cubic form\footnote{It is particularly remarkable that in case of octonions the definition
of such Norm is invariant w.r.t changes in the multiplication order
of the the term $\left(x_{1}x_{2}\right)x_{3}+\overline{x}_{3}\left(\overline{x}_{2}\overline{x}_{1}\right)$.} 
\begin{equation}
N\left(X\right)=\lambda_{1}\lambda_{2}\lambda_{3}-\lambda_{1}x_{1}\overline{x}_{1}-\lambda_{2}x_{2}\overline{x}_{2}-\lambda_{3}x_{3}\overline{x}_{3}+\left(x_{1}x_{2}\right)x_{3}+\overline{x}_{3}\left(\overline{x}_{2}\overline{x}_{1}\right).\label{eq:Cubic Norm}
\end{equation}
Then, the cubic norm is Jordan admissible and yields to the Jordan
product 
\begin{equation}
X\circ Y=\frac{1}{2}\left(XY+YX\right),
\end{equation}
with trace given by $\text{Tr}\left(X\right)=\lambda_{1}+\lambda_{2}+\lambda_{3}$,
and bilinear-trace form given by $\left(X,Y\right)=\text{Tr}\left(X\circ Y\right)$
and where the juxtaposition denotes the ordinary matrix product throughout.
Finally, the sharp map of the generic element $X$ is of the form
\begin{equation}
X^{\#}=\left(\begin{array}{ccc}
\lambda_{2}\lambda_{3}-n\left(x_{1}\right) & \overline{x}_{2}\overline{x}_{1}-\lambda_{3}x_{3} & x_{3}x_{1}-\lambda_{2}\overline{x}_{2}\\
x_{1}x_{2}-\lambda_{3}\overline{x}_{3} & \lambda_{1}\lambda_{3}-n\left(x_{2}\right) & \overline{x}_{3}\overline{x}_{2}-\lambda_{1}x_{1}\\
\overline{x}_{1}\overline{x}_{3}-\lambda_{2}x_{2} & x_{2}x_{3}-\lambda_{1}\overline{x}_{1} & \lambda_{1}\lambda_{2}-n\left(x_{3}\right)
\end{array}\right).\label{eq:Sharp-map-X}
\end{equation}

The definition of the sharp map $\#$ allows a classification of elements
of the Jordan algebra through the notion of\emph{ rank }of the element.
This notion is due to Jacobson \cite{Jac} and will be useful in characterizing
Veronese vectors in following chapters. Given a cubic norm $N$, an
element $X$ is said to be 
\begin{itemize}
\item of rank $3$, i.e. $rank\left(X\right)=3$ if and only if $N\left(X\right)\neq0$;
\item of rank $2$, i.e. $rank\left(X\right)=2$ if and only if $N\left(X\right)=0$
and $X^{\#}\neq0$; 
\item of rank $1$, i.e. $rank\left(X\right)=1$ if and only if $X\neq0$
and $X^{\#}=0$; 
\item finally $X=0$ is of rank 0.
\end{itemize}

\subsection{Albert algebras }

We now turn our attention to the Albert algebra, which belongs to
the category of exceptional Jordan algebras. Over the field of real
numbers $\mathbb{R}$ these algebras amount to four inequivalent exceptional
algebras, i.e. $\mathfrak{J}_{3}\left(\mathbb{O}\right)$, $\mathfrak{J}_{3}\left(\mathbb{O}_{s}\right)$,
$\mathfrak{J}_{2,1}\left(\mathbb{O}\right)$ and $\mathfrak{J}_{2,1}\left(\mathbb{O}_{s}\right)$,
that we will treat in a common setup following \cite{Jac60}. Let
$\mathfrak{J}$ be the 27-dimensional real vector space $H_{3}\left(\mathbb{O};\gamma_{1},\gamma_{2},\gamma_{3}\right)$
of matrices of the form 
\begin{equation}
X=\left(\begin{array}{ccc}
\lambda_{1} & x_{3} & \gamma_{1}^{-1}\gamma_{3}\overline{x}_{2}\\
\gamma_{2}^{-1}\gamma_{1}\overline{x}_{3} & \lambda_{2} & x_{1}\\
x_{2} & \gamma_{3}^{-1}\gamma_{2}\overline{x}_{1} & \lambda_{3}
\end{array}\right),
\end{equation}
 with $\gamma_{1},\gamma_{2},\gamma_{3}\in\left\{ \pm1\right\} $,
$\lambda_{1},\lambda_{2},\lambda_{3}\in\mathbb{R}$, $x_{1},x_{2},x_{3}\in\mathbb{O}$,
endowed with the following bilinear product 
\begin{equation}
X\circ Y=\frac{1}{2}\left(X\eta Y+Y\eta X\right),
\end{equation}
for every $X,Y\in\mathfrak{J}$, where $\eta=\text{diag}\left(\gamma_{1},\gamma_{2},\gamma_{3}\right)$
and the juxtaposition is the ordinary matrix product. The algebra
$\mathfrak{J}$ is a Jordan algebra with cubic norm which is the octonionic
version of the determinant (3.13) of $X$ 
\begin{equation}
N\left(X\right)=\lambda_{1}\lambda_{2}\lambda_{3}-\gamma_{3}^{-1}\gamma_{2}\lambda_{1}n\left(x_{1}\right)-\gamma_{1}^{-1}\gamma_{3}\lambda_{2}n\left(x_{2}\right)-\gamma_{2}^{-1}\gamma_{1}\lambda_{3}n\left(x_{3}\right)+\left\langle x_{1}x_{2},\overline{x}_{3}\right\rangle ,\label{eq:Gen Cubic norm}
\end{equation}
where the definition of octonionic norm has been used (the octonionic
product is abbreviated by juxtaposition). Applying previous constructions
we obtain 

\begin{equation}
\text{Tr\ensuremath{\left(X\right)}}=\gamma_{1}\lambda_{1}+\gamma_{2}\lambda_{2}+\gamma_{3}\lambda_{3},\label{eq:TracciaTensoriale}
\end{equation}
and the\emph{ Freudenthal product}, which is twice the polarisation
of the sharp map $\#$, i.e. $x\times y=\left(x+y\right)^{\#}-x^{\#}-y^{\#}$,
is given by 
\begin{equation}
X\times Y=\frac{1}{2}\left(2X\circ Y-X\text{Tr\ensuremath{\left(Y\right)}}-Y\text{Tr\ensuremath{\left(X\right)}}+\left(\text{Tr\ensuremath{\left(X\right)}}\text{Tr\ensuremath{\left(Y\right)}}-\text{\text{Tr\ensuremath{\left(X\circ Y\right)}}}\right)\text{\textbf{1}}\right).
\end{equation}

We then have the following notable relations involving the Jordan
product and Freudenthal product
\begin{align}
X\circ Y & =Y\circ X,\\
X\times Y & =Y\times X,\\
2\left(X\times X\right)\circ X & =N\left(X\right)\text{\textbf{1}},\label{eq:CubicNorm from Freud}\\
4\left(X\times X\right)\times\left(X\times X\right) & =N\left(X\right)X,
\end{align}
 where $X,Y\in\mathfrak{J}$ and $\boldsymbol{1}$ is the three-by-three
identity matrix.

For a broader understanding of the Albert algebra in contexts beyond
Hurwitz algebras, especially where the concept of Hermitian matrices
isn't straightforwardly extended, it is insightful to observe how
the Jordan and Freudenthal products operate on basis elements. Using
this, and through (\ref{eq:Cubic Norm}), we can define a cubic norm
and a sharp map from $X^{\#}=\frac{1}{2}\left(X\times X\right)$. 

Then, let us consider the following elements of $\mathfrak{J}$, i.e.
\begin{align}
e_{1} & =\left(\begin{array}{ccc}
1 & 0 & 0\\
0 & 0 & 0\\
0 & 0 & 0
\end{array}\right),e_{2}=\left(\begin{array}{ccc}
0 & 0 & 0\\
0 & 1 & 0\\
0 & 0 & 0
\end{array}\right),e_{3}=\left(\begin{array}{ccc}
0 & 0 & 0\\
0 & 0 & 0\\
0 & 0 & 1
\end{array}\right),\\
\iota_{12}\left(a\right) & =\left(\begin{array}{ccc}
0 & a & 0\\
\gamma_{2}^{-1}\gamma_{1}\overline{a} & 0 & 0\\
0 & 0 & 0
\end{array}\right),\iota_{13}\left(a\right)=\left(\begin{array}{ccc}
0 & 0 & \gamma_{1}^{-1}\gamma_{3}\overline{a}\\
0 & 0 & 0\\
a & 0 & 0
\end{array}\right),\iota_{23}\left(a\right)=\left(\begin{array}{ccc}
0 & 0 & 0\\
0 & 0 & a\\
0 & \gamma_{3}^{-1}\gamma_{2}\overline{a} & 0
\end{array}\right),\label{eq:iij Elements base}
\end{align}
so that we have the following decomposition
\begin{equation}
\mathfrak{J}=\mathbb{R}e_{1}\mathbb{\oplus}\mathbb{R}e_{2}\mathbb{\oplus}\mathbb{R}e_{3}\mathbb{\oplus}\mathfrak{J}_{12}\mathbb{\oplus}\mathfrak{J}_{23}\mathbb{\oplus}\mathfrak{J}_{13},
\end{equation}
for $\mathfrak{J}_{ij}=\left\{ \iota_{ij}\left(a\right)\in\mathfrak{J}\right\} $.
Then we have that the Jordan product is defined by $e_{i}\circ e_{j}=\delta_{ij}e_{i},$
and the following relations \cite{Jac60}
\begin{equation}
\begin{cases}
e_{i}\circ\iota_{ij}\left(a\right)=\frac{1}{2}\iota_{ij}\left(a\right), & ,\,\,\,\\
\iota_{ij}\left(a\right)\circ\iota_{ij}\left(a\right)=\gamma_{j}^{-1}\gamma_{i}n\left(a\right)\left(e_{i}+e_{j}\right),\\
\iota_{ij}\left(a\right)\circ\iota_{jk}\left(b\right)=\iota_{ik}\left(ab\right),
\end{cases}
\end{equation}
 for $a,b\in\mathbb{O}$. 

It is now important to state that in this framework, the four inequivalent
Jordan algebras arise, i.e. $\mathfrak{J}_{3}\left(\mathbb{O}\right)$,
$\mathfrak{J}_{3}\left(\mathbb{O}_{s}\right)$, $\mathfrak{J}_{2,1}\left(\mathbb{O}\right)$
and $\mathfrak{J}_{2,1}\left(\mathbb{O}_{s}\right)$. Those four exceptional
Jordan algebras, which are the only exceptional Jordan algebras over
the field of reals are given by $H_{3}\left(\mathbb{O};1,1,1\right)$,
$H_{3}\left(\mathbb{O}_{s};1,1,1\right)$, $H_{3}\left(\mathbb{O};1,1,-1\right)$
and $H_{3}\left(\mathbb{O}_{s};1,1,-1\right)$ respectively. 

\section{\label{sec:real-forms-of}Real forms of Lie algebras }

A profound relationship exists between octonionic algebras and exceptional
Lie groups. Indeed, we have already recalled how the group of automorphisms
of the octonions is the compact form of $G_{2}$, while the one of
the split-octonions is its non-compact form $G_{2(2)}$. Historically,
was Cartan in \cite{Car14} to consider the automorphisms of the algebra
of octonions $\text{Aut}\left(\mathbb{O}\right)$ as a model of $\text{G}_{2}$
and thus the algebra $\mathfrak{g}_{2}$ as its derivations algebra,
i.e. $\mathfrak{der}\left(\mathbb{O}\right)$. However, a systematic
exploration of this connection only commenced when Chevalley and Schafer
\cite{ChSch} proved that the exceptional algebra $\mathfrak{f}_{4}$
could be obtained as the algebra of derivation of the Albert algebra,
i.e., $\mathfrak{der}\left(\mathfrak{J}_{3}\left(\mathbb{O}\right)\right)$.
Subsequent research\cite{Fr54,Tits,Freud 1965,Vinberg} culminated
in the Tits-Freudenthal construction. The Tits-Freudenthal magic square
provides a systematic methodology to deduce all exceptional Lie algebras
from Jordan algebras and tensor products over octonions. Alternatively,
similar outcomes can be achieved through different methods. For instance,
using Kantor triple systems and their generalizations, all exceptional
Lie algebras can be derived from algebraic constructions over octonions
\cite{Pa06,Pa08}. 

To briefly summarize relations between the Albert algebra $\mathfrak{J}_{3}\left(\mathbb{O}\right)$
and exceptional Lie algebras, we have the the following suggestive
set of results

\begin{equation}
\begin{array}{cc}
\mathfrak{der}\left(\mathbb{O}\right) & =\mathfrak{g}_{2},\\
\mathfrak{der}\left(\mathfrak{J}_{3}\left(\mathbb{O}\right)\right) & =\mathfrak{f}_{4},\\
\mathfrak{str}_{0}\left(\mathfrak{J}_{3}\left(\mathbb{O}\right)\right) & =\mathfrak{e}_{6},\\
\mathfrak{conf}\left(\mathfrak{J}_{3}\left(\mathbb{O}\right)\right) & =\mathfrak{e}_{7},\\
\mathfrak{qconf}\left(\mathfrak{J}_{3}\left(\mathbb{O}\right)\right) & =\mathfrak{e}_{8},
\end{array}
\end{equation}
where $\mathfrak{der}\left(A\right)$, $\mathfrak{str}_{0}\left(A\right)$,
$\mathfrak{conf}\left(A\right)$ and $\mathfrak{qconf}\left(A\right)$
are the derivation algebra, the reduced structure algebra, the conformal
algebra and the quasiconformal algebra of $A$ respectively (for a
definition of conformal and quasiconformal algebra see \cite{GKN}).

As for the exceptional Lie groups, it is worth noting that, at the
present days, the Albert algebra is considered as the main tool in
defining exceptional Lie groups such as $\text{F}_{4}$ and $\text{E}_{6}$
(e.g. see \cite{Yokota}). In the next section we will see that the
octonionic planes allow also an equivalent but purely geometric definition
of such groups. Since proofs and calculations are not significantly
altered, for the remain of this section we will consider $\gamma_{1}=\gamma_{2}=\gamma_{3}=1$
for the algebra of octonions $\mathbb{O}$ and only state results
for the other relevant cases.

\subsection{Lie groups of type $\text{F}_{4}$ }

Let us define the compact form of $\text{F}_{4}$ as the automorphism
group of the Albert algebra $\mathfrak{J}_{3}\left(\mathbb{O}\right)$,
here simply denoted as $\mathfrak{J}$, i.e.

\begin{equation}
\text{F}_{4}=\left\{ \varphi\in\text{End}\left(\mathfrak{J}\right):\varphi\left(X\circ Y\right)=\varphi\left(X\right)\circ\varphi\left(Y\right)\text{ and }X,Y\in\mathfrak{J}\right\} .\label{eq:Definition F4}
\end{equation}
We then show that it is a Lie group with simple Lie algebra of dimension
$52$ and thus consistent with notation used for Cartan's classification
of semisimple Lie algebras. 

To show that the group of automorphisms of the Albert algebra $\mathfrak{J}_{3}\left(\mathbb{O}\right)$
is a Lie group we proceed showing that any automorphism is trace-preserving,
i.e. $\text{Tr}\left(\varphi\left(X\right)\right)=\text{Tr}\left(X\right)$.
Rewriting (\ref{eq:CubicNorm from Freud}) as 
\begin{equation}
\left(X\circ X\right)\circ X-\text{Tr}\left(X\right)X^{2}+\frac{1}{2}\left(\text{Tr}\left(X\right)^{2}-\text{Tr}\left(X^{2}\right)\right)X=N\left(X\right)\text{\textbf{1}},
\end{equation}
applying it for $\varphi\left(X\right)$ and then applying $\varphi^{-1}$
to the obtained expression we have
\begin{equation}
\left(X\circ X\right)\circ X-\text{Tr}\left(\varphi\left(X\right)\right)X^{2}+\frac{1}{2}\left(\text{Tr}\left(\varphi\left(X\right)\right)^{2}-\text{Tr}\left(\varphi\left(X\right)^{2}\right)\right)X=N\left(\varphi\left(X\right)\right)\text{\textbf{1}},
\end{equation}
 that once is subtracted to the previous gives
\begin{equation}
\left(\text{Tr}\left(\varphi\left(X\right)\right)-\text{Tr}\left(X\right)\right)X^{2}+\frac{1}{2}\left(\text{Tr}\left(X\right)^{2}-\text{Tr}\left(\varphi\left(X\right)\right)^{2}+\text{Tr}\left(\varphi\left(X\right)^{2}\right)-\text{Tr}\left(X^{2}\right)\right)=\left(N\left(\varphi\left(X\right)\right)-N\left(X\right)\right)\text{\textbf{1}}.
\end{equation}
Setting $X$as elements of the basis in (\ref{eq:iij Elements base})
we obtain that once $X=\varphi\left(\iota_{ij}\left(a\right)\right)$
then $\text{Tr}\left(\varphi\left(\iota_{ij}\left(a\right)\right)\right)=0$
as $\text{Tr}\left(\iota_{ij}\left(a\right)\right)=0$, and $\text{Tr}\left(\varphi\left(e_{i}\right)\right)=\text{Tr}\left(e_{i}\right)$
so that in fact we obtain 
\begin{equation}
\text{Tr}\left(\varphi\left(X\right)\right)=\text{Tr}\left(X\right),\label{eq:Tr(phi(X))=00003DTr(X)}
\end{equation}
for every $\varphi\in\text{F}_{4}$ and $X\in\mathfrak{J}$. 

Since (\ref{eq:Tr(phi(X))=00003DTr(X)}) and since $S\left(X,Y\right)=\text{Tr\ensuremath{\left(X\circ Y\right)}}$
then we have that 
\begin{equation}
S\left(\varphi\left(X\right),\varphi\left(Y\right)\right)=S\left(X,Y\right),
\end{equation}
which can be shown to be an equivalent condition with $\varphi\left(X\circ Y\right)=\varphi\left(X\right)\circ\varphi\left(Y\right)$
(e.g. \cite{Yokota}). This means that $\text{F}_{4}$ can be characterized
as a subgroup of the special orthogonal group of $\mathbb{R}^{27}$,
i.e. 
\begin{equation}
\text{SO}\left(27\right)=\left\{ \varphi\in\text{End}\left(\mathfrak{J}\right):S\left(\varphi\left(X\right),\varphi\left(Y\right)\right)=S\left(X,Y\right)\right\} ,
\end{equation}
and therefore is a compact Lie group. 

The real dimension of $\text{F}_{4}$ is achieved considering its
Lie algebra $\mathfrak{f}_{4}$, i.e. the derivation algebras of the
Albert algebra $\mathfrak{J}$ \cite{ChSch}. The mapping of the form
\begin{equation}
\text{ad}\left(X\right):Y\longrightarrow XY-YX,
\end{equation}
where $X$ is a traceless skew-Hermitian matrix, i.e. $\overline{X}^{t}=-X$
and $\text{Tr}\left(X\right)=0$, are derivations that form a subspace
of the Lie algebra $\mathfrak{der}\left(\mathfrak{J}\right)$ usually
denoted by $\mathfrak{sa}\left(3\right)$. Then, any derivation $D$
of the octonionic algebra $\mathbb{O}$ generates a derivation of
$\mathfrak{J}$ just appling the derivation $D$ to each entry of
the traceless skew-Hermitian matrix $X$. Thus, we have that the decomposition
of $\mathfrak{J}$ in 
\begin{equation}
\mathfrak{der}\left(\mathfrak{J}_{3}\left(\mathbb{O}\right)\right)=\mathfrak{der}\left(\mathbb{O}\right)\oplus\mathfrak{sa}\left(3\right)=\mathfrak{g}_{2}\oplus\mathfrak{sa}\left(3\right).
\end{equation}
Here the elements of $\mathfrak{sa}\left(3\right)$ are of the form
\begin{equation}
A=\left(\begin{array}{ccc}
a_{1}^{1} & a_{2}^{1} & -\overline{a}_{3}^{1}\\
-\overline{a}_{2}^{1} & a_{2}^{2} & a_{3}^{2}\\
a_{3}^{1} & -\overline{a}_{3}^{2} & a_{3}^{3}
\end{array}\right)
\end{equation}
with $a_{3}^{3}=-\left(a_{1}^{1}+a_{2}^{2}\right)$ and $\text{Re}\left(a_{1}^{1}\right)=\text{Re}\left(a_{2}^{2}\right)=0$.
We therefore have 3 coefficient of dimension $8$, $2$ entries of
dimension $7$ and therefore $\text{dim}_{\mathbb{R}}\mathfrak{sa}\left(3\right)=38$
so that 
\begin{equation}
\mathfrak{\text{dim}_{\mathbb{R}}\mathfrak{der}\left(\mathfrak{J}_{3}\left(\mathbb{O}\right)\right)}\cong52=38+14.
\end{equation}

Some subgroups of $\text{F}_{4}$ are relevant for the next chapter.
Obviously in $\text{F}_{4}$ we can identify a subgroup isomorphic
to $\text{G}_{2}$ which is given by the automorphisms of the octonions
$\mathbb{O}$, so that for every automorphism $\varphi\in\text{G}_{2}$
we have an induced automorphism $\widetilde{\varphi}$ in $\text{F}_{4}$
defined as 
\begin{equation}
\widetilde{\varphi}\left(e_{i}\right)=e_{i},\,\,\,\,\widetilde{\varphi}\left(\iota_{ij}\left(a\right)\right)=\iota_{ij}\left(\varphi\left(a\right)\right),
\end{equation}
for every $a\in\mathbb{O}$. We thus have that $\text{F}_{4}$ contains
a (non-maximal - with commutant $\text{SU}(2)$- , and non-symmetric)
subgroup isomorphic to $\text{G}_{2}$, i.e. $\text{F}_{4}\supset\text{G}_{2}$.
Moreover, relaxing previous conditions and imposing just that the
automorphisms fix the elements $e_{1},e_{2}$ and $e_{3}$, i.e.
\begin{equation}
\varphi\left(e_{i}\right)=e_{i},\text{\,\,\,\,\,}i=1,2,3,
\end{equation}
we obtain a subgroup that is isomorphic to $\text{Spin}\left(8\right)$,
i.e. 
\begin{equation}
\left(\text{F}_{4}\right)_{\text{diag}}=\left\{ \varphi\in\text{F}_{4}:\varphi\left(e_{i}\right)=e_{i},i=1,2,3\right\} \cong\text{Spin}\left(8\right).
\end{equation}
 Finally, imposing just $\varphi\left(e_{1}\right)=e_{1}$ we obtain
a subgroup that is isomorphic to $\text{Spin}\left(9\right)$, i.e. 

\begin{align}
\left(\text{F}_{4}\right)_{e_{1}} & =\left\{ \varphi\in\text{F}_{4}:\varphi\left(e_{1}\right)=e_{1}\right\} \cong\text{Spin}\left(9\right).\label{eq:(F4)e1=00003DSpin(9)}
\end{align}
Complete proofs of previous statements can be found in \cite{Yo68}.
In conclusion, for $\text{F}_{4}$ we have the following chain of
maximal subgroups inclusions
\begin{equation}
\begin{array}{ccccccccccccc}
\text{F}_{4} & \longrightarrow & \text{Spin}\left(9\right) & \longrightarrow & \text{Spin}\left(8\right) & \longrightarrow & \text{Spin}\left(7\right) & \longrightarrow & \text{G}_{2(-14)} & \longrightarrow & \text{SU}\left(3\right) & \longrightarrow & \text{SU}\left(2\right)\end{array}.\label{eq:F4SubGroup}
\end{equation}

For our purposes, also real forms of $\text{F}_{4}$ are relevant.
All that was developed for the automorphisms of the Albert algebra
$\mathfrak{J}_{3}\left(\mathbb{O}\right)$, yielding to the compact
form of $\text{F}_{4}$, sometimes denoted as $\text{F}_{4\left(-52\right)}$,
can be reconsidered for the automorphism group of $\mathfrak{J}_{2,1}\left(\mathbb{O}\right)$
and $\mathfrak{J}_{3}\left(\mathbb{O}_{s}\right)$ with similar results.
All resulting groups are Lie groups of the $\text{F}_{4}$ type but
with different signature in the Cartan-Killing form \cite{Jac60}.
More specifically we have the following identifications

\begin{align}
\text{F}_{4\left(-52\right)} & =\left\{ \varphi\in\text{End}\left(\mathfrak{J}_{3}\left(\mathbb{O}\right)\right):\varphi\left(X\circ Y\right)=\varphi\left(X\right)\circ\varphi\left(Y\right)\right\} ,\label{eq:Definition F4-1}\\
\text{F}_{4\left(-20\right)} & =\left\{ \varphi\in\text{End}\left(\mathfrak{J}_{2,1}\left(\mathbb{O}\right)\right):\varphi\left(X\circ Y\right)=\varphi\left(X\right)\circ\varphi\left(Y\right)\right\} ,\\
\text{F}_{4\left(4\right)} & =\left\{ \varphi\in\text{End}\left(\mathfrak{J}_{3}\left(\mathbb{O}_{s}\right)\right):\varphi\left(X\circ Y\right)=\varphi\left(X\right)\circ\varphi\left(Y\right)\right\} \\
 & =\left\{ \varphi\in\text{End}\left(\mathfrak{J}_{2,1}\left(\mathbb{O}_{s}\right)\right):\varphi\left(X\circ Y\right)=\varphi\left(X\right)\circ\varphi\left(Y\right)\right\} .\nonumber 
\end{align}
 While the chain of subgroups in (\ref{eq:F4SubGroup}) is referred
to the real compact form $\text{F}_{4(-52)}$, the minimally non-compact
real form of $\text{F}_{4}$ resulting from the Lorentzian Jordan
algebra $\mathfrak{J}_{2,1}\left(\mathbb{O}\right)$, i.e. $\text{F}_{4\left(-20\right)}$,
enjoys the following chain of maximal inclusions
\begin{equation}
\begin{array}{ccccccccccc}
\text{F}_{4\left(-20\right)} & \longrightarrow & \text{Spin}\left(9\right) & \longrightarrow & \text{Spin}\left(8\right) & \longrightarrow & \text{Spin}\left(7\right) & \longrightarrow & \text{G}_{2(-14)} & \longrightarrow & \text{SU}\left(3\right)\\
 & \searrow &  & \nearrow &  & \nearrow\\
 &  & \text{Spin}\left(8,1\right) & \longrightarrow & \text{Spin}\left(7,1\right) & \longrightarrow & \text{Spin}\left(6,1\right)
\end{array}
\end{equation}
On the other hand $\mathfrak{J}_{3}\left(\mathbb{O}_{s}\right)$ and
$\mathfrak{J}_{2,1}\left(\mathbb{O}_{s}\right)$ yields to the same
real form of $\text{F}_{4}$, i.e. the (maximally non-compact) split
form $\text{F}_{4(4)}$, for which we have the chain of maximal inclusions
\begin{equation}
\begin{array}{ccccccccccc}
 &  &  &  &  &  & \text{Spin}\left(5,2\right) &  &  &  & \text{SU}\left(2,1\right)\\
 &  &  &  &  & \nearrow &  &  &  & \nearrow\\
\text{F}_{4\left(4\right)} & \longrightarrow & \text{Spin}\left(5,4\right) & \longrightarrow & \text{Spin}\left(5,3\right) & \longrightarrow & \text{Spin}\left(4,3\right) & \longrightarrow & \text{G}_{2(2)}\\
 &  &  & \searrow &  & \nearrow &  &  &  & \searrow\\
 &  &  &  & \text{Spin}\left(4,4\right) &  &  &  &  &  & \text{SL}\left(3,\mathbb{R}\right)
\end{array}
\end{equation}

\subsection{Lie groups of type $\text{E}_{6}$ }

If the exceptional groups of type $\text{F}_{4}$ are realised as
the automorphism group of the Albert algebra $\mathfrak{J}_{3}\left(\mathbb{O}\right)$,
or equivalently as the trace-preserving transformations of $\mathfrak{J}_{3}\left(\mathbb{O}\right)$,
and the exceptional groups of type $\text{E}_{6}$ are obtained as
the symmetry groups of the cubic norm, i.e. as the determinant preserving,
transformations of the Albert algebra $\mathfrak{J}_{3}\left(\mathbb{O}\right)$.
Indeed, the complex form of $\text{E}_{6}$ is defined as 
\begin{equation}
\text{E}_{6}^{\mathbb{C}}=\left\{ \varphi\in\text{End}\left(\mathfrak{J}_{3}^{\mathbb{C}}\left(\mathbb{O}\right)\right):N\left(\varphi\left(X\right)\right)=N\left(X\right)\text{ }X\in\mathfrak{J}_{3}^{\mathbb{C}}\left(\mathbb{O}\right)\right\} ,\label{eq:Complex E6}
\end{equation}
where $N\left(X\right)$ is the cubic norm defined in (\ref{eq:Cubic Norm})
which is the determinant of the element $X$ in its matrix realisation
(\ref{eq:Hermitian element}). To obtain the minimally non-compact
real form of $\text{E}_{6}$, one needs either $\mathfrak{J}_{3}\left(\mathbb{O}\right)$
or $\mathfrak{J}_{2,1}\left(\mathbb{O}\right)$ (see \cite{CCM}),
i.e.,
\begin{align}
\text{E}_{6\left(-26\right)} & =\left\{ \varphi\in\text{End}\left(\mathfrak{J}_{3}\left(\mathbb{O}\right)\right):N\left(\varphi\left(X\right)\right)=N\left(X\right)\right\} \label{eq:Definition E6}\\
 & =\left\{ \varphi\in\text{End}\left(\mathfrak{J}_{2,1}\left(\mathbb{O}\right)\right):N\left(\varphi\left(X\right)\right)=N\left(X\right)\right\} .\nonumber 
\end{align}
 which is thus seen as the group of cubic-norm-preserving linear transformations
of $\mathfrak{J}_{3}\left(\mathbb{O}\right)$. The dimension of such
group and, thus, its classification is the result of the theorem by
Chevalley-Schafer that in 1950 triggered the research for links between
Jordan algebras and Lie algebras, i.e. 
\begin{thm}
\emph{(Chevalley-Schafer \cite{ChSch})} The exceptional simple Lie
algebra $\mathfrak{e}_{6}$ is given by 
\begin{equation}
\mathfrak{e}_{6}=\mathfrak{der}\left(\mathfrak{J}_{3}\left(\mathbb{O}\right)\right)\oplus\mathfrak{J}_{3}^{0}\left(\mathbb{O}\right),
\end{equation}
where $\mathfrak{J}_{3}^{0}\left(\mathbb{O}\right)$ identifies the
traceless elements of the Albert algebra.
\end{thm}

As a corollary of the previous, we have that 
\begin{align}
\text{dim}\mathfrak{e}_{6} & =\text{dim}\left(\mathfrak{der}\left(\mathfrak{J}_{3}\left(\mathbb{O}\right)\right)\right)\oplus\text{dim}\left(\mathfrak{J}_{3}^{0}\left(\mathbb{O}\right)\right)\\
 & =52+26=78.\nonumber 
\end{align}
The other real form of $\text{E}_{6}$ of interest for our work is
the split form $\text{E}_{6\left(6\right)}$, obtainable by the means
of either $\mathfrak{J}_{3}\left(\mathbb{O}\right)$ or $\mathfrak{J}_{2,1}\left(\mathbb{O}\right)$,
i.e., 

\begin{align}
\text{E}_{6\left(6\right)} & =\left\{ \varphi\in\text{End}\left(\mathfrak{J}_{3}\left(\mathbb{O}_{s}\right)\right):N\left(\varphi\left(X\right)\right)=N\left(X\right)\right\} \\
 & =\left\{ \varphi\in\text{End}\left(\mathfrak{J}_{2,1}\left(\mathbb{O}_{s}\right)\right):N\left(\varphi\left(X\right)\right)=N\left(X\right)\right\} .\nonumber 
\end{align}

As for their relations with real forms of $\text{F}_{4}$ the respective
split forms are related by a maximal group embedding,
\begin{equation}
\text{E}_{6\left(6\right)}\supset\text{F}_{4\left(4\right)}.
\end{equation}
On the other hand, the real form $\text{E}_{6\left(-26\right)}$ contains
both $\text{F}_{4\left(-20\right)}$ and $\text{F}_{4\left(-52\right)}$,
i.e. 

\begin{equation}
\text{E}_{6\left(-26\right)}\supset\text{F}_{4\left(-20\right)}\text{ and }\text{E}_{6\left(-26\right)}\supset\text{F}_{4\left(-52\right)},
\end{equation}
the latter being its maximal compact subgroup.

\section{\label{sec:Veronese-formulation-of}Veronese formulation of the octonionic
planes }

It is a common practice defining a projective plane over an associative
division algebra $\mathbb{K}\in\left\{ \mathbb{R},\mathbb{C},\mathbb{H}\right\} $
from the tridimensional vector space over the given algebra, e.g.
$\mathbb{K}^{3}\setminus\left\{ 0\right\} $, and then define the
projective plane as the quotient 
\begin{equation}
\mathbb{K}P^{2}=\left(\mathbb{K}^{3}\setminus\left\{ 0\right\} \right)/\sim,
\end{equation}
where $x\sim y$ if $x$ and $y$ are non-zero vectors multiple through
the scalar field, i.e. $\lambda x=y$, $\lambda\in\mathbb{K}$, $\lambda\neq0$,
$x,y\in\mathbb{K}^{3}\setminus\left\{ 0\right\} $. However, this
definition encounters complications for the octonionic case, i.e.,
$\mathbb{O}P^{2}$, due to the non-associative nature of the algebra
of octonions $\mathbb{O}$. Specifically, we have 
\begin{equation}
\left(\lambda\mu\right)x\neq\lambda\left(\mu x\right)
\end{equation}
when $\lambda,\mu\in\mathbb{O}$. So, if we attempt to define the
equivalence relation in the previous manner, $x\sim y=\lambda x$
and $y\sim z=\mu y$ do not imply $x\sim z$ since
\begin{equation}
z=\mu\left(\lambda x\right)\neq\left(\lambda\mu\right)x.
\end{equation}
Consequently, the aforementioned relation does not qualify as an equivalence
relation, making the quotient ill-defined.

Although it became apparent that the traditional definition was unsuitable
for the octonionic projective plane, at the beginning of the XX century
it was already understood that such a projective plane existed. Through
investigations carried out by Cartan \cite{Car15}, Jordan, Wigner
\& von Neumann \cite{Jordan} and Freudenthal \cite{Fr54} on the
geometry of the projective plane, four equivalent definitions for
the projective plane over $\mathbb{K}$ emerged: 1) as a completion
of the affine plane $\mathbb{K}A^{2}$; 2) as the trace-1 idempotent
elements of the rank-three Jordan algebra\footnote{when K is over the reals, one might also consider $\mathfrak{J}_{2,1}\left(\mathbb{K}\right)$.}
$\mathfrak{J}_{3}\left(\mathbb{K}\right)$; 3) as a coset manifold
obtained with a specific isometry and isotropy group and, finally,
4) through some \emph{ad hoc} but self-contained definitions. 

For our purposes, a powerful self-contained definition is the one
that rely on Veronese coordinates. This direct construction was pioneered
by Salzmann et al. in \cite{Compact Projective}, and we will show
in Sec. \ref{sec:Realisations-through-symmetric} is easily generalisable
to symmetric composition algebras. Remarkably, the following diagramis
commutative:
\[
\begin{array}{ccc}
\text{Veronese} & \Longleftrightarrow & \text{Compl. Affine}\\
\Updownarrow & \mathbb{\mathbb{K}}P^{2} & \Updownarrow\\
\text{Coset Manif.} & \Longleftrightarrow & \underset{\text{rank-1 idempotents}}{\text{Jordan}\,\mathfrak{J}_{3}\left(\mathbb{\mathbb{K}}\right)}
\end{array}
\]
These constructions are robust and can be adapted to all Hurwitz algebras,
provided we cautiously acknowledge that the term \textquotedblleft projective\textquotedblright{}
should be reinterpreted in the context of non-division (i.e., split)
Hurwitz algebras, as the axioms of projective geometry are not satisfied.

Now, we will see how three of them apply to the octonionic and split-octonionic
planes: the definition through the aid of Veronese vectors; the algebraic
definition through the Albert algebras; the definition as homogenous
spaces. On the other hand, we leave the first approach, i.e. the projective
plane as a completion of the affine plane $\mathbb{K}A^{2}$, for
future publications since that it would constitute an unnecessary
deviation from the topic we wish to develop.

\section{\label{sec:The-projective-octonionic}The projective octonionic plane }

An \emph{incidence plane} $P^{2}$ is given by the triple $\left\{ \mathscr{P},\mathscr{L},\mathscr{R}\right\} $
where $\mathscr{P}$ is the set of points of the plane, $\mathscr{L}$
the set of lines and $\mathscr{R}$ are the incidence relations of
poins and lines. For $P^{2}$ to be projective, $\mathscr{R}$ must
satisfy the following properties:
\begin{enumerate}
\item Any two distinct points are incident to a unique line.
\item Any two distinct lines are incident with a unique point. 
\item (\emph{non-degenera}cy) There exist four points such that no three
are incident one another.
\end{enumerate}
Then, points passing through the same line are said \emph{collinear},
and lines passing through the same point are said \emph{concurrent}.
An orderd set of three distinct non-collinear points is said a \emph{triangle},
and if the set is of four point no three of which are collinear, then
it is called a \emph{quadrangle}.\emph{ }

We now define the octonionic projective plane $\mathbb{O}P^{2}$.
Let $V\cong\mathbb{O}^{3}\times\mathbb{R}^{3}$ be a real vector space,
with elements of the form 
\[
\left(x_{\nu};\lambda_{\nu}\right)_{\nu}\coloneqq\left(x_{1},x_{2},x_{3};\lambda_{1},\lambda_{2},\lambda_{3}\right)
\]
where $x_{\nu}\in\mathbb{O}$, $\lambda_{\nu}\in\mathbb{R}$ and $\nu=1,2,3$.
A vector $w\in V$ is called \emph{Veronese} \cite{Compact Projective}
if 

\begin{align}
\lambda_{1}\overline{x}_{1} & =x_{2}\cdot x_{3},\,\,\lambda_{2}\overline{x}_{2}=x_{3}\cdot x_{1},\,\,\lambda_{3}\overline{x}_{3}=x_{1}\cdot x_{2}\label{eq:Ver-1}\\
n\left(x_{1}\right) & =\lambda_{2}\lambda_{3},\,n\left(x_{2}\right)=\lambda_{3}\lambda_{1},n\left(x_{3}\right)=\lambda_{1}\lambda_{2}.\label{eq:Ver-2}
\end{align}
Let $H\subset V$ be the subset of Veronese vectors. If $w=\left(x_{\nu};\lambda_{\nu}\right)_{\nu}$
is a Veronese vector then also its real multiple $\mu w=\mu\left(x_{\nu};\lambda_{\nu}\right)_{\nu}$
is a Veronese vector, that is all real multiples of $w$ are Veronese,
i.e. $\mathbb{R}w\subset H$. We then define \emph{points} $\mathscr{P}_{\mathbb{O}}$
of the octonionic projective plane $\mathbb{O}P^{2}$ as the 1-dimensional
subspace $\mathbb{R}w$ where $w$ is a Veronese vector, i.e.
\begin{equation}
\mathscr{P}_{\mathbb{O}}=\left\{ \mathbb{R}w:w\in H\smallsetminus\left\{ 0\right\} \right\} .\label{eq:Octonionic point}
\end{equation}

\begin{rem}
A point in the projective plane is defined as the equivalence class
$\mathbb{R}w$ of the Veronese vector $w$, but, in order to make
an explicit relation between points in the projective plane and rank-1
idempotent elements of the Albert algebra, we will choose as representative
the vector $v=\left(y_{\nu};\xi_{\nu}\right)_{\nu}\in\mathbb{R}w$
such that $\xi_{1}+\xi_{2}+\xi_{3}=1$.

We then define the set of \emph{projective lines} $\mathscr{L}_{\mathbb{O}}$
of $\mathbb{O}P^{2}$ as the vectors orthogonal to the points $\mathbb{R}w$.
In order to define the lines of the incidence plane, we need to define
a bilinear form that allows to evaluate orthogonality. Let $\beta$
be the symmetric bilinear form over $\mathbb{O}^{3}\times\mathbb{R}^{3}$
defined as
\end{rem}

\begin{equation}
\beta\left(w_{1},w_{2}\right)=\stackrel[\nu=1]{3}{\sum}\left(\left\langle x_{\nu}^{1},x_{\nu}^{2}\right\rangle +\lambda_{\nu}^{1}\lambda_{\nu}^{2}\right)
\end{equation}
where $w_{1}=\left(x_{\nu}^{1};\lambda_{\nu}^{1}\right)_{\nu}$,$w_{2}=\left(x_{\nu}^{2};\lambda_{\nu}^{2}\right)_{\nu}\in\mathbb{O}^{3}\times\mathbb{R}^{3}$
and $\left\langle x,y\right\rangle =\overline{x}\cdot y+\overline{y}\cdot x$.
Then, for every point $\mathbb{R}w$ in $\mathscr{P}_{\mathbb{O}}$,
corresponding to the Veronese vector $w$, one defines a line $\ell$
in $\mathbb{O}P^{2}$ as the orthogonal space through $\beta$, i.e.,
\begin{equation}
\ell\coloneqq w^{\perp}=\left\{ z\in\mathbb{O}^{3}\times\mathbb{R}^{3}:\beta\left(z,w\right)=0\right\} ,\label{eq:Octonionic lines}
\end{equation}
so that $\mathscr{L}_{\mathbb{O}}$ is given by the set of lines defined
in (\ref{eq:Octonionic lines}). Explicitly, $\beta\left(w_{1},w_{2}\right)=0$
if and only if
\begin{equation}
\overline{x}_{1}^{1}x_{1}^{2}+\overline{x}_{1}^{2}x_{1}^{1}+\overline{x}_{2}^{1}x_{2}^{2}+\overline{x}_{2}^{2}x_{2}^{1}+\overline{x}_{3}^{1}x_{3}^{2}+\overline{x}_{3}^{2}x_{3}^{1}+\lambda_{1}^{1}\lambda_{1}^{2}+\lambda_{2}^{1}\lambda_{2}^{2}+\lambda_{3}^{1}\lambda_{3}^{2}=0.\label{eq:Explicit-line}
\end{equation}

Finally, the incidence relations $\mathscr{R}$ are given by $p\in\mathscr{P}_{\mathbb{O}}$
is incident to $\ell\in\mathscr{L}_{\mathbb{O}}$ iff $p$ belongs
to $\ell$, i.e. $p\in\ell$, so that the incidence relation is actually
given by the inclusion $\subseteq$. Thus, the \emph{octonionic projective
plane} $\mathbb{O}P^{2}$ is defined as the incidence plane resulting
from the triple
\begin{equation}
\mathbb{O}P^{2}\coloneqq\left\{ \mathscr{P}_{\mathbb{O}},\mathscr{L}_{\mathbb{O}},\subseteq\right\} .
\end{equation}

It is worth noting that, by its definition, the bilinear form $\beta$
implicitly defines a polarity on the octonionic projective plane $\mathbb{O}P^{2}$,
i.e. an order two bijective map between lines and points. Indeed,
the map $\pi^{+}$, called the \emph{elliptic polarity}, that corresponds
point to lines and lines to points, i.e.
\begin{equation}
\pi^{+}\left(w\right)=w^{\perp},\pi^{+}\left(w^{\perp}\right)=w,
\end{equation}
 is defined by the bilinear form $\beta\left(\cdot,\cdot\right)$,
so that
\begin{align}
\pi^{+}: & w\longrightarrow\left\{ \beta\left(\cdot,w\right)=0\right\} ,\label{eq:elliptic polarity}\\
 & \ell\longrightarrow w,
\end{align}
where the line $\ell$ is given by the set $\left\{ \beta\left(\cdot,w\right)=0\right\} $.
We thus have that any statement related to projective points it corresponds
to a statement related to projective lines, as the classical principle
of duality between points and lines states, i.e.
\begin{thm}
\noun{(}Principle of projective duality\noun{)} In any non-degenerate
projective plane,statements that are valid for points of the planes
are dually valid for its lines.
\end{thm}

It is also worth noting that in this specific case, since we are leading
with a composition division algebra, Veronese conditions can be reduced
to a minimal set with the following
\begin{prop}
\label{prop:17Dimensional}In case of $V\cong\mathbb{O}^{3}\times\mathbb{R}^{3}$,
conditions (\ref{eq:Ver-1}) and (\ref{eq:Ver-2}) are equivalent
to the following indipendent conditions 
\end{prop}

\begin{align}
\lambda_{1}\overline{x}_{1} & =x_{2}\cdot x_{3},\label{eq:Indip Ver}\\
n\left(x_{2}\right) & =\lambda_{3}\lambda_{1},n\left(x_{3}\right)=\lambda_{1}\lambda_{2}.\label{eq:Indip Ver2}
\end{align}

\begin{proof}
Since $\mathbb{O}$ are a division algebra, let multiply (\ref{eq:Indip Ver})
on the left by $\overline{x}_{2}$ in order to obtain 
\begin{equation}
\lambda_{1}\overline{x}_{2}\cdot\overline{x}_{1}=n\left(x_{2}\right)x_{3},
\end{equation}
then applying (\ref{eq:Indip Ver2}) we obtain 
\begin{equation}
\lambda_{1}\overline{x}_{2}\cdot\overline{x}_{1}=\lambda_{3}\lambda_{1}x_{3},
\end{equation}
which by octonionic conjugation gives back $\lambda_{3}\overline{x_{3}}=x_{1}x_{2}$.
Similarly, multiplying on the right by $\overline{x}_{3}$ we obtain
$\lambda_{2}\overline{x_{2}}=x_{3}x_{1}$. Finally, since $\mathbb{O}$
is composition, then 
\begin{equation}
n\left(x_{1}\right)=x_{1}\cdot\overline{x}_{1}=\lambda_{1}^{-2}n\left(\overline{x}_{2}\cdot\overline{x}_{3}\right)=\lambda_{2}\lambda_{3},
\end{equation}
which is the last Veronese condition.
\end{proof}
As anticipated above, with only minor modifications to the previous
arguments, one can define other octonionic planes using Veronese-like
coordinates. Indeed, let $\widetilde{V}=\mathbb{K}^{3}\times\mathbb{R}^{3}$
where $\mathbb{K}$ can either be the octonions or the split-octonions,
i.e. $\mathbb{K}\in\left\{ \mathbb{O},\mathbb{O}_{s}\right\} $. For
a triple $\left(\gamma_{1},\gamma_{2},\gamma_{3}\right)$ with $\gamma_{1},\gamma_{2},\gamma_{3}\in\left\{ \pm1\right\} $,
we define a Veronese-like vector as a non-zero vector $\left(x_{\nu};\lambda_{\nu}\right)_{\nu}\in\widetilde{V}$
such that 
\begin{align}
\lambda_{1}\overline{x}_{1} & =\gamma_{3}\gamma_{2}^{-1}x_{2}\cdot x_{3},\,\,\lambda_{2}\overline{x}_{2}=\gamma_{1}\gamma_{3}^{-1}x_{3}\cdot x_{1},\,\,\,\lambda_{3}\overline{x}_{3}=\gamma_{2}\gamma_{1}^{-1}x_{1}\cdot x_{2},\label{eq:VerGen1-1}\\
n\left(x_{1}\right) & =\gamma_{3}\gamma_{2}^{-1}\lambda_{2}\lambda_{3},\,\,n\left(x_{2}\right)=\gamma_{1}\gamma_{3}^{-1}\lambda_{3}\lambda_{1},\,\,\,n\left(x_{3}\right)=\gamma_{2}\gamma_{1}^{-1}\lambda_{1}\lambda_{2},\label{eq:VerGen2-1}
\end{align}
where $x_{\nu}\in\mathbb{K}$ and $\lambda_{\nu}\in\mathbb{R}$ for
$\nu=1,2,3$. We call such subspace $\widetilde{H}$. Notably, one
observe that if $x_{\nu}\in\mathbb{O}_{s}$ is a zero divisor, then
$n_{s}\left(x_{\nu}\right)=0$ so that in case of $x_{2}\cdot x_{3}=0$,
the Veronese conditions imply $\lambda_{1}=0$, since $n\left(x_{2}\right)=\lambda_{3}\lambda_{1}=0$
and $n\left(x_{3}\right)=\lambda_{1}\lambda_{2}=0$ and the vector
is non-zero. \emph{Points} $\mathscr{P}$ of the octonionic incidence
plane are defined as the 1-dimensional subspaces $\mathbb{R}w$ where
$w$ is a Veronese-like vector, i.e.
\begin{equation}
\mathscr{P}=\left\{ \mathbb{R}w:w\in\widetilde{H}\smallsetminus\left\{ 0\right\} \right\} .\label{eq:Octonionic point-1}
\end{equation}
 As for the line we proceed defining the symmetric bilinear form $\beta$
as
\begin{equation}
\beta\left(x,y\right)=\stackrel[\nu=1]{3}{\sum}\left(\gamma_{\nu+2}^{-1}\gamma_{\nu+1}\left\langle x_{\nu},y_{\nu}\right\rangle +\lambda_{\nu}\mu_{\nu}\right),
\end{equation}
where $x=\left(x_{\nu};\lambda_{\nu}\right)_{\nu}$ and $y=\left(y_{\nu};\mu_{\nu}\right)_{\nu}$
and indices are taken modulo $3$. Then, the set of the lines $\mathscr{L}$
are given as the orthogonal space through $\beta$ to a Veronese vector
$x$, i.e.

\begin{equation}
\ell\coloneqq x^{\perp}=\left\{ y\in\mathbb{O}^{3}\times\mathbb{R}^{3}:\beta\left(y,x\right)=0\right\} ,
\end{equation}
which explicitly happens when

\begin{equation}
\gamma_{3}^{-1}\gamma_{2}\left\langle x_{1},y_{1}\right\rangle +\gamma_{1}^{-1}\gamma_{3}\left\langle x_{2},y_{2}\right\rangle +\gamma_{2}^{-1}\gamma_{1}\left\langle x_{3},y_{3}\right\rangle +\lambda_{1}\mu_{1}+\lambda_{2}\mu_{2}+\lambda_{3}\mu_{3}=0.\label{eq:Explicit-line-1}
\end{equation}
Lastly, the incidence relations remain unchanged and are always determined
by inclusion.

It can be easily shown (see \cite[Sec. 16]{Compact Projective}) that
for $\left(\gamma_{1},\gamma_{2},\gamma_{3}\right)=\left(1,1,1\right)$
and $\mathbb{K}=\mathbb{O}$ we obtain a definition of the octonionic
projective plane $\mathbb{O}P^{2}$. On the other hand, setting $\left(\gamma_{1},\gamma_{2},\gamma_{3}\right)=\left(1,1,-1\right)$
but leaving $\mathbb{K}=\mathbb{O}$, we obtain the \emph{octonionic
hyperbolic plane} $\mathbb{O}H^{2}$. On the other hand, switching
to split-octonions we obtain the \emph{split-octonionic projective
}or \emph{hyperbolic plane}, i.e $\mathbb{O}_{s}P^{2}$ or $\mathbb{O}_{s}H^{2}$,
depending on $\gamma_{3}$ being $1$ or $-1$, respectively.
\begin{table}
\centering{}%
\begin{tabular}{|c|c|c|c|}
\hline 
Plane & $\mathbb{K}$ & $\left(\gamma_{1},\gamma_{2},\gamma_{3}\right)$ & $\mathfrak{J}$\tabularnewline
\hline 
\hline 
$\mathbb{O}P^{2}$ & $\mathbb{O}$ & $\left(1,1,1\right)$ & $\mathfrak{J}_{3}\left(\mathbb{O}\right)$\tabularnewline
\hline 
$\mathbb{O}H^{2}$ & $\mathbb{O}$ & $\left(1,1,-1\right)$ & $\mathfrak{J}_{2,1}\left(\mathbb{O}\right)$\tabularnewline
\hline 
$\mathbb{O}_{s}P^{2}$ & $\mathbb{O}_{s}$ & $\left(1,1,1\right)$ & $\mathfrak{J}_{3}\left(\mathbb{O}_{s}\right)$\tabularnewline
\hline 
$\mathbb{O}_{s}H^{2}$ & $\mathbb{O}_{s}$ & $\left(1,1,-1\right)$ & $\mathfrak{J}_{2,1}\left(\mathbb{O}_{s}\right)$\tabularnewline
\hline 
\end{tabular}\caption{\label{tab:Octonionic-incidence-planes,}Octonionic incidence planes,
algebra of definition, the triple $\left(\gamma_{1},\gamma_{2},\gamma_{3}\right)$
and the corresponding Albert algebra $\mathfrak{J}$.}
\end{table}

\section{\label{sec:Collineations-and-real}Collineations and real forms of
Lie algebras }

The algebraic representation of the octonionic incidence planes is
established by considering the mapping $\psi$ between $\widetilde{V}$
and the Albert agebra $\mathfrak{J}$ given by
\begin{equation}
\psi:\left(x_{\nu};\lambda_{\nu}\right)_{\nu}\mapsto\left(\begin{array}{ccc}
\lambda_{1} & x_{3} & \gamma_{1}^{-1}\gamma_{3}\overline{x}_{2}\\
\gamma_{2}^{-1}\gamma_{1}\overline{x}_{3} & \lambda_{2} & x_{1}\\
x_{2} & \gamma_{3}^{-1}\gamma_{2}\overline{x}_{1} & \lambda_{3}
\end{array}\right),\label{eq:mappingVer-Jordan}
\end{equation}
where $\mathfrak{J}$ can be $\mathfrak{J}_{3}\left(\mathbb{O}\right)$,
$\mathfrak{J}_{2,1}\left(\mathbb{O}\right)$,$\mathfrak{J}_{3}\left(\mathbb{O}_{s}\right)$
and $\mathfrak{J}_{2,1}\left(\mathbb{O}\right)$ depending on the
plane considered (see Table \ref{tab:Octonionic-incidence-planes,}).
Then, since the cubic form $N\left(X\right)$, as defined in (\ref{eq:Gen Cubic norm}),
and the $\#$ map of a generic element is given by 

\begin{equation}
X^{\#}=\left(\begin{array}{ccc}
\lambda_{2}\lambda_{3}-\gamma_{3}^{-1}\gamma_{2}n\left(x_{1}\right) & \gamma_{1}^{-1}\gamma_{2}\overline{x}_{2}\overline{x}_{1}-\lambda_{3}x_{3} & x_{3}x_{1}-\gamma_{1}^{-1}\gamma_{3}\lambda_{2}\overline{x}_{2}\\
x_{1}x_{2}-\gamma_{2}^{-1}\gamma_{1}\lambda_{3}\overline{x}_{3} & \lambda_{1}\lambda_{3}-\gamma_{1}^{-1}\gamma_{3}n\left(x_{2}\right) & \gamma_{2}^{-1}\gamma_{3}\overline{x}_{3}\overline{x}_{2}-\lambda_{1}x_{1}\\
\gamma_{3}^{-1}\gamma_{1}\overline{x}_{1}\overline{x}_{3}-\lambda_{2}x_{2} & x_{2}x_{3}-\gamma_{3}^{-1}\gamma_{2}\lambda_{1}\overline{x}_{1} & \lambda_{1}\lambda_{2}-\gamma_{2}^{-1}\gamma_{1}n\left(x_{3}\right)
\end{array}\right),
\end{equation}
we have the following characterization 
\begin{align}
\mathbb{O}P^{2} & \cong\left\{ X\in\mathfrak{J}_{3}\left(\mathbb{O}\right):X^{\#}=0,\text{Tr}\left(X\right)=1\right\} ,\\
\mathbb{O}H^{2} & \cong\left\{ X\in\mathfrak{J}_{2,1}\left(\mathbb{O}\right):X^{\#}=0,\text{Tr}\left(X\right)=1\right\} ,\\
\mathbb{O}_{s}P^{2} & \cong\left\{ X\in\mathfrak{J}_{3}\left(\mathbb{O}_{s}\right):X^{\#}=0,\text{Tr}\left(X\right)=1\right\} ,\\
\mathbb{O}_{s}H^{2} & \cong\left\{ X\in\mathfrak{J}_{2,1}\left(\mathbb{O}_{s}\right):X^{\#}=0,\text{Tr}\left(X\right)=1\right\} ,
\end{align}
where $\text{Tr}\left(X\right)=\gamma_{1}\lambda_{1}+\gamma_{2}\lambda_{2}+\gamma_{3}\lambda_{3}$
as defined in (\ref{eq:TracciaTensoriale}).
\begin{rem}
Note that alternative definitions are possible, e.g., for references
see \cite{Sp60,SpVe63,Manogue}.
\end{rem}

As one might expect, in all cases the collineation group is a Lie
group of type $\text{E}_{6}$ while the isometry group is given by
a Lie group of type $\text{F}_{4}$. From the algebraic equivalent
of the octonionic planes, from Albert algebras one derives \cite[p.105]{Yokota}
that the automorphism groups are in fact real forms of $\text{E}_{6}$
while the isometry groups are given by the real forms of $\text{F}_{4}$
(see Table \ref{tab:Relations-between-incidence}).
\begin{rem}
It is interesting to notice that the non-compact real forms $E_{6(2)}$
and $E_{6(-14)}$ do not appear as collineation groups of some real
form of the octonionic projective plane. This is due to the fact that
they are not directly related to real forms of the Albert algebras
(See e.g. \cite[Sec. 3.10-11]{Yokota}), or equivalently to the fact
that they do not appear in the second row of the relevant real form
of the Freudenthal-Tits Magic Square (see Tables 12 and 13 of \cite{CCM}).
\end{rem}

\begin{table}
\begin{centering}
\begin{tabular}{|c|c|c|c|}
\hline 
Incidence plane $\mathscr{\mathcal{P}}^{2}\left(\mathbb{K}\right)$ & $\text{Aut}\left(\mathscr{\mathcal{P}}^{2}\left(\mathbb{K}\right)\right)$ & $\text{Iso}\left(\mathscr{\mathcal{P}}^{2}\left(\mathbb{K}\right)\right)$ & $\Gamma\diamondsuit\left(\mathscr{\mathcal{P}}^{2}\left(\mathbb{K}\right)\right)$\tabularnewline
\hline 
\hline 
$\mathbb{O}P_{\mathbb{C}}^{2}$ & $\text{E}_{6}^{\mathbb{C}}$ & $\text{F}_{4}^{\mathbb{C}}$ & $\text{G}_{2}^{\mathbb{C}}$\tabularnewline
\hline 
$\mathbb{O}P^{2}$ & $\text{E}_{6\left(-26\right)}$ & $\text{F}_{4\left(-52\right)}$ & $\text{G}_{2\left(-14\right)}$\tabularnewline
\hline 
$\mathbb{O}_{s}P^{2}$ & $\text{E}_{6\left(6\right)}$ & $\text{F}_{4\left(4\right)}$ & $\text{G}_{2\left(2\right)}$\tabularnewline
\hline 
$\mathbb{O}_{s}H^{2}$ & $\text{E}_{6\left(6\right)}$ & $\text{F}_{4\left(4\right)}$ & $\text{G}_{2\left(2\right)}$\tabularnewline
\hline 
$\mathbb{O}H^{2}$ & $\text{E}_{6\left(-26\right)}$ & $\text{F}_{4\left(-20\right)}$ & $\text{G}_{2\left(-14\right)}$\tabularnewline
\hline 
\end{tabular}
\par\end{centering}
\caption{\label{tab:Relations-between-incidence}Relations between incidence
plane, the collineation group $\text{Aut}\left(\mathscr{\mathcal{P}}^{2}\left(\mathbb{K}\right)\right)$,
the isometry group $\text{Iso}\left(\mathscr{\mathcal{P}}^{2}\left(\mathbb{K}\right)\right)$
and the stabilizer of a non degenerate quadrangle $\Gamma\diamondsuit\left(\mathscr{\mathcal{P}}^{2}\left(\mathbb{K}\right)\right)$
for $\mathbb{K}\in\left\{ \mathbb{O},\mathbb{O}_{s}\right\} $ and
$\mathscr{\mathcal{P}}\in\left\{ P,H\right\} $. See \cite{Compact Projective,Sa17,Sa08}
for references.}

\end{table}

\begin{center}
\begin{figure}
\centering{}\includegraphics[scale=0.5]{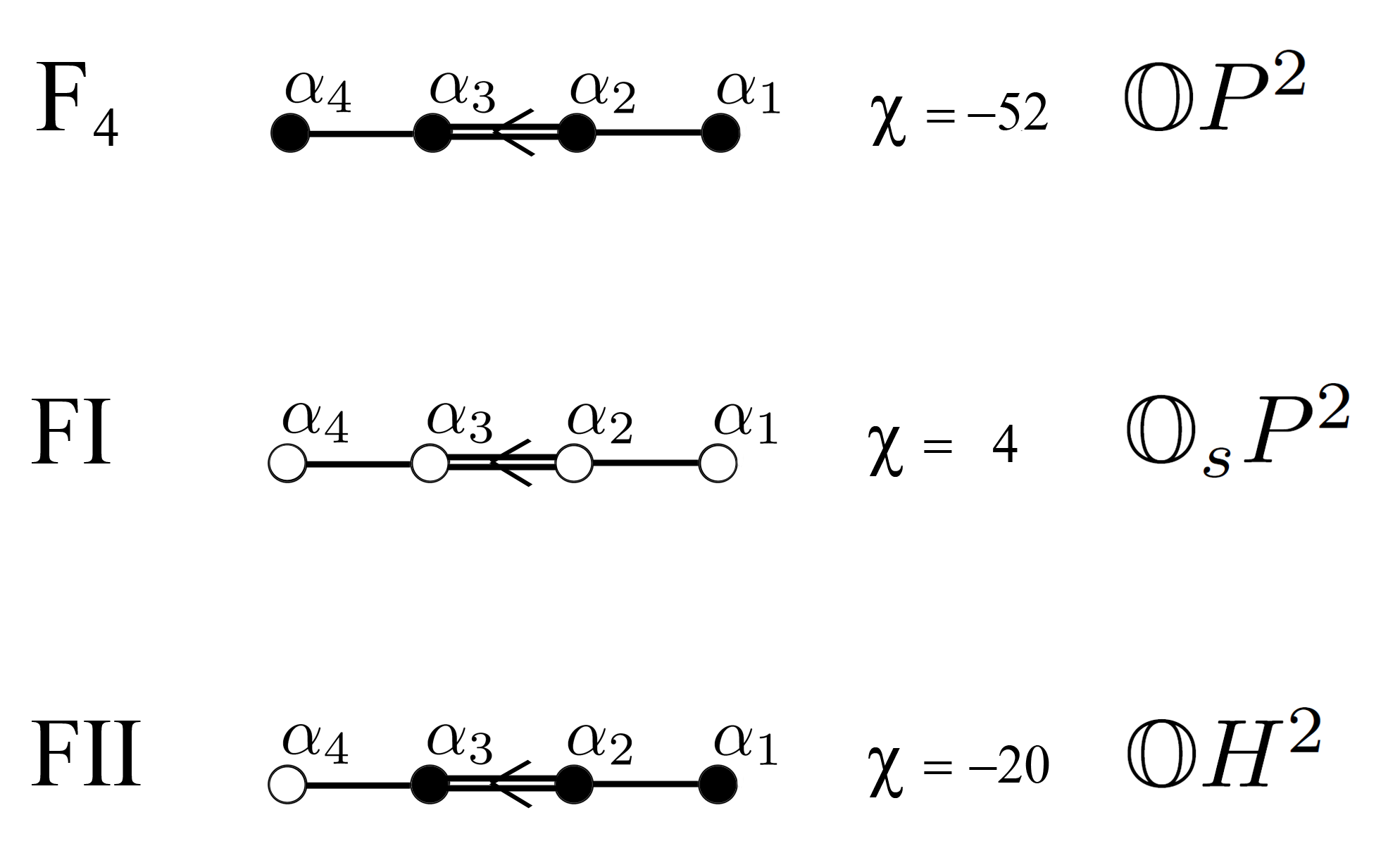}\caption{Satake diagrams of real forms of $\text{F}_{4}$, their character
\textgreek{q} and corresponding plane of which they are the isometry
group.}
\end{figure}
\par\end{center}

Finally, inspired by the table in the previous section and by Rosenfeld's
work \cite{Rosenfeld-1993,CCMA-Magic} we can classify al possible
octonionic planes considering all possible coset manifolds arising
from real forms of $\text{F}_{4}$ considered as isometry and with
real dimension $16$. Starting from the complexification of the Cayley
plane (for an explicit construction see Sec. \ref{sec:Complex-Cayley-plane})
\begin{eqnarray}
\mathbb{O}P_{\mathbb{C}}^{2} & \cong & \frac{\text{F}_{4}^{\mathbb{C}}}{\text{Spin}_{\mathbb{C}}\left(9\right)},
\end{eqnarray}
we define three different real forms of the plane: a totally compact
real coset manifold, that in fact we have seen being isomorphic to
$\mathbb{O}P^{2}$; a totally non-compact which can be shown to be
isomorphic to the hyperbolic octonionic plane $\mathbb{O}H^{2}$.
Those octonionic planes will be defined taking as isometry group $\text{F}_{4\left(-52\right)}$
and $\text{F}_{4\left(-20\right)}$ while the last real form $\text{F}_{4\left(4\right)}$,
since it contains as subgroup $\text{G}_{2\left(2\right)}$ and not
$\text{G}_{2\left(-14\right)}$, will yield to projective planes on
the split-octonionic algebra $\mathbb{O}_{s}$, i.e.

\begin{align}
\mathbb{O}P^{2} & \cong\frac{\text{F}_{4(-52)}}{\text{Spin}\left(9\right)},\\
\mathbb{O}H^{2} & \cong\frac{\text{F}_{4(-20)}}{\text{Spin}\left(9\right)},\\
\mathbb{O}_{s}P^{2} & \cong\mathbb{O}_{s}H^{2}\cong\frac{\text{F}_{4(4)}}{\text{Spin}\left(5,4\right)}.
\end{align}

The classification of all the planes proceeds through the use of type
and character, where the type identifies the cardinality of non-compact
and compact generators, i.e. $\left(\#_{nc},\#_{c}\right)$, and the
character $\chi$ is given by the difference between the two, i.e.
$\chi=\#_{nc}-\#_{c}$. Using this classification we then see that,
out of the three different octonionic planes, $\mathbb{O}P^{2}$ is
totally compact of type $(0,16)$ and has character $\chi=-16$ (the
classical \emph{Cayley-Moufang plane}) ; $\mathbb{O}H^{2}$ is totally\footnote{It is here worth remarking that a pseudo-Riemannian form of the hyperbolic
plane $\mathbb{O}H^{2}$ also exists, given by the symmetric coset
$\text{F}_{4(-20)}$/Spin(8,1), which is again of type $(8,8)$ and
thus of character $\chi=0$. Seemingly, such a coset does not have
an immediate interpretation as a \emph{locus} inside a real form of
the Albert algebra.} non-compact of type $\left(16,0\right)$ and has character $\chi=-16$
(the\emph{ hyperbolic octonionic plane}) ; and $\mathbb{O}_{s}P^{2}$
is isomorphic to $\mathbb{O}_{s}H^{2}$ and is of type $\left(8,8\right)$
and has character $\chi=0$ . 

\section{\label{sec:Realisations-through-symmetric}Realisations through symmetric
composition algebras}

In this section, we introduce two novel realizations of both the hyperbolic
octonionic plane $\mathbb{O}H^{2}$ and of the split-octonionic projective
plane $\mathbb{O}_{s}P^{2}$. These realisations are build on our
prior work \cite{CMZ}, were we constructed an independent realisation
of the Cayley-Moufang plane through the use of symmetric composition
and division algebras such as the Okubo algebra $\mathcal{O}$ and
the paraoctonionic algebra $p\mathbb{O}$. In the cited work, we adapted
the Veronese conditions (\ref{eq:Ver-1}) and (\ref{eq:Ver-2}) to
suit the definition of projective planes over symmetric composition
and division algebras, specifically, the Okubic projective plane $\mathcal{O}P^{2}$
and the paraoctonionic projective plane $p\mathbb{O}P^{2}$. We demonstrated
that these newly defined projective planes are isomorphic to the projective
completion of the corresponding affine planes $\overline{\mathscr{A}^{2}}\left(\mathcal{O}\right)\cong\mathcal{O}P^{2}$
and $\overline{\mathscr{A}^{2}}\left(p\mathbb{O}\right)\cong p\mathbb{O}P^{2}$.
Most strikingly, we established that even though both Okubo and paraoctonionic
algebras are not alternative, the resulting projective planes are
nevertheless both isomorphic to the octonionic projective plane, that
is,
\begin{equation}
\mathcal{O}P^{2}\cong p\mathbb{O}P^{2}\cong\mathbb{O}P^{2}.
\end{equation}
This result is particularly noteworthy because, as shown in Table
\ref{tab:Synoptic-table-of}, the three algebras, while deeply interconnected,
are not isomorphic. Specifically, the octonions $\mathbb{O}$ are
alternative and unital, paraoctonions $p\mathbb{O}$ are neither alternative
nor unital but they do have a para-unit, while the Okubo algebra $\mathcal{O}$
is non-alternative, non-unital, and it only has idempotents elements.
It is also worth highlighting that the Okubo algebra $\mathcal{O}$
is the one with smaller automorphism group among these algebras (for
a summary of the property of these algebras see Table \ref{tab:Synoptic-table-of}).
\begin{table}
\centering{}%
\begin{tabular}{|c|c|c|c|}
\hline 
Property & $\mathbb{O}$ & $p\mathbb{O}$ & $\mathcal{O}$\tabularnewline
\hline 
\hline 
Unital & Yes & No & No\tabularnewline
\hline 
Paraunital & Yes & Yes & No\tabularnewline
\hline 
Alternative & Yes & No & No\tabularnewline
\hline 
Flexible & Yes & Yes & Yes\tabularnewline
\hline 
Composition & Yes & Yes & Yes\tabularnewline
\hline 
Division & Yes & Yes & Yes\tabularnewline
\hline 
Automorphism & $\text{G}_{2}$ & $\text{G}_{2}$ & $\text{SU}\left(3\right)$\tabularnewline
\hline 
\end{tabular}\caption{\label{tab:Synoptic-table-of}Synoptic table of the algebraic properties
of octonions $\mathbb{O}$, paraoctonions $p\mathbb{O}$ and the real
Okubo algebra $\mathcal{O}$.}
\end{table}

In this section, we build upon our earlier work to present new realizations
of the hyperbolic octonionic plane $\mathbb{O}H^{2}$ and of the split-octonionic
projective plane $\mathbb{O}_{s}P^{2}$ through non-unital and non-alternative
algebras. We believe that these novel realizations are both intriguing
and potentially useful, particularly in the context of Particle Physics,
as detailed in \cite{CMZ-2}. The structure of this section is as
follows. First, we provide a concise overview of symmetric composition
algebras and explain how to derive the Okubo algebra $\mathcal{O}$
and the split-Okubo algebra $\mathcal{O}_{s}$ along with the paraoctonionic
$p\mathbb{O}$ and split-paraoctonionic algebra $p\mathbb{O}_{s}$
(see \cite{ElDuque Comp} for a more comprehensive treatment of all
these algebras). Then, we present the general setting for hyperbolic
and projective planes as in (\ref{eq:VerGen1-1}) and (\ref{eq:VerGen2-1})
but suitably modified for symmetric composition algebras. Finally,
we present the isomorphism between the Okubic, paraoctonionic and
octonionic hyperbolic planes, i.e. $\mathcal{O}H^{2}\cong p\mathbb{O}H^{2}\cong\mathbb{O}H^{2}$,
along with the isomorphism between the split-Okubo, split-paraoctonionic
and split-octonionic projective planes, i.e. $\mathcal{O}_{s}P^{2}\cong p\mathbb{O}_{s}P^{2}\cong\mathbb{O}_{s}P^{2}.$

\subsection{Non-unital composition algebras}

As briefly reviewed in Section \ref{sec:octonions-and-split-octonions},
normed algebras endowed with a norm that respect the multiplicative
structure, i.e., $n\left(x\cdot y\right)=n\left(x\right)n\left(y\right)$,
are called composition algebras. Composition algebras are divided
in \emph{unital}, i.e. where a unit element $1$ exists such that
$x\cdot1=1\cdot x=x$, \emph{para-unital}, i.e. where an involution
$x\longrightarrow\overline{x}$ and an element called paraunit $\boldsymbol{1}$
exists such that $x\cdot\boldsymbol{1}=\boldsymbol{1}\cdot x=\overline{x}$,
and, finally, \emph{non-unital}, i.e. referring here to algebras that
do not possess neither a unit nor a paraunit element. 
\begin{table}
\begin{centering}
\begin{tabular}{|c|c|c|c|c|c|c|c|c|c|c|c|c|}
\hline 
\textbf{Hurwitz} & \textbf{O.} & \textbf{C.} & \textbf{A.} & \textbf{Alt.} & \textbf{F.} &  & \textbf{p-Hurwitz} & \textbf{O.} & \textbf{C.} & \textbf{A.} & \textbf{Alt.} & \textbf{F.}\tabularnewline
\hline 
\hline 
$\mathbb{R}$ & Yes & Yes & Yes & Yes & Yes &  & $p\mathbb{R}\cong\mathbb{R}$ & Yes & Yes & Yes & Yes & Yes\tabularnewline
\hline 
$\mathbb{C}$, $\mathbb{C}_{s}$ & No & Yes & Yes & Yes & Yes &  & $p\mathbb{C}$, $p\mathbb{C}_{s}$ & No & Yes & No & No & Yes\tabularnewline
\hline 
$\mathbb{H}$,$\mathbb{H}_{s}$ & No & No & Yes & Yes & Yes &  & $p\mathbb{H}$,$p\mathbb{H}_{s}$ & No & No & No & No & Yes\tabularnewline
\hline 
$\mathbb{O}$,$\mathbb{O}_{s}$ & No & No & No & Yes & Yes &  & $p\mathbb{O}$,$p\mathbb{O}_{s}$ & No & No & No & No & Yes\tabularnewline
\hline 
\end{tabular}
\par\end{centering}
\centering{}{\small{}\bigskip{}
}\caption{\emph{\label{tab:Hurwitz-para-Hurwitz-1}On the left,} we have summarized
the algebraic properties, i.e. totally ordered (O), commutative (C),
associative (A), alternative (Alt), flexible (F), of all Hurwitz algebras,
namely $\mathbb{R},\mathbb{C},\mathbb{H}$ and $\mathbb{O}$ along
with their split counterparts $\mathbb{C}_{s},\mathbb{H}_{s},\mathbb{O}_{s}$.
\emph{On the right}, we have summarized the algebraic properties of
all para-Hurwitz algebras, namely $p\mathbb{R},p\mathbb{C},p\mathbb{H}$
and $p\mathbb{O}$ accompanied by their split counterparts $p\mathbb{C}_{s},p\mathbb{H}_{s},p\mathbb{O}_{s}$.}
\end{table}

In fact, as a consequence of the Generalized Hurwitz Theorem, only
sixteen composition algebras exist for any given field (see \cite{ElDuque Comp,ZSSS}
for reference): seven are unital and called \emph{Hurwitz algebras},
comprising the division algebras $\mathbb{R},\mathbb{C},\mathbb{H},\mathbb{O}$
along with their split companions $\mathbb{C}_{s},\mathbb{H}_{s},\mathbb{O}_{s}$;
other seven are paraunital, closely related to the Hurwitz algebras
and termed \emph{para-Hurwitz algebras} (see Table \ref{tab:Hurwitz-para-Hurwitz-1});
finally, there are two composition algebras, one division and one
split, that are both non-unital, non-alternative, and 8-dimensional,
known as the \emph{Okubo algebras} $\mathcal{O}$ and $\mathcal{O}_{s}$.
In the next subsection we will be interested in the four para-unital
or non-unital composition algebras of dimension eight: the para-Hurwitz
algebras of paraoctonions $p\mathbb{O}$ and split-paraoctonions $p\mathbb{O}_{s}$,
along with the Okubo and split-Okubo algebras $\mathcal{O}$ and $\mathcal{O}_{s}$. 

\subsection{\label{subsec:Conjugation-and-the}The algebras of paraoctonions
and split-paraoctonions}

In unital composition algebras there exists a canonical involution,
namely an order-two anti-homomorphism, known as \emph{conjugation}.
This can be defined using the orthogonal projection of the unit element
as
\begin{equation}
x\mapsto\overline{x}=\left\langle x,1\right\rangle 1-x.\label{eq:coniugazione}
\end{equation}
This canonical involution has the distinctive property of being an
antihomomorphism with respect to the product, i.e., $\overline{x\cdot y}=\overline{y}\cdot\overline{x},$
and the basic property of $x\cdot\overline{x}=n\left(x\right)1$. 

Given the order-two antihomomorphism of the conjugation over the Hurwitz
algebra of octonions $\left(\mathbb{O},\cdot,n\right)$ we can now
obtain a para-Hurwitz algebra defining a new product
\begin{equation}
x\bullet y=\overline{x}\cdot\overline{y},
\end{equation}
 for every $x,y\in\mathbb{O}$. The new algebra $\left(\mathbb{O},\bullet,n\right)$
is again a composition algebra, in fact the para-Hurwitz algebra,
called of paraoctonions, denoted with $p\mathbb{O}$. As expected,
the algebra does not have a unit but only a para-unit, i.e. $\boldsymbol{1}\in p\mathbb{O}$
such that $\boldsymbol{1}\bullet x=x\bullet\boldsymbol{1}=\overline{x}.$
It is worth noting that also the algebra of paraoctonions is a division
algebra, since if 
\begin{equation}
x\bullet y=\overline{x}\cdot\overline{y}=0,
\end{equation}
then either $\overline{x}$ or $\overline{y}$ are zero, and thus
either $x$ or $y$ are zero. 

Proceeding in a similar manner but with the split-octonions $\mathbb{O}_{s}$,
we obtain a new algebra $\left(\mathbb{O}_{s},\bullet,n\right)$ that,
again, is a composition algebra of dimension eight. This algebra is
non-alternative, and, as the para-octonionic algebra, is a para-unital
algebra, but (as the split-octonions) it is not a division algebra
and it is thus called the split-paraoctonionic algebra $p\mathbb{O}_{s}$.

\subsection{The Okubo algebra}

Given by the following map 
\begin{equation}
\begin{array}{cc}
\tau\left(\text{e}_{k}\right) & =\text{e}_{k},k=0,1,3,7\\
\tau\left(\text{e}_{2}\right) & =-\frac{1}{2}\left(\text{e}_{2}-\sqrt{3}\text{e}_{5}\right),\\
\tau\left(\text{e}_{5}\right) & =-\frac{1}{2}\left(\text{e}_{5}+\sqrt{3}\text{e}_{2}\right),\\
\tau\left(\text{e}_{4}\right) & =-\frac{1}{2}\left(\text{e}_{4}-\sqrt{3}\text{e}_{6}\right),\\
\tau\left(\text{e}_{6}\right) & =-\frac{1}{2}\left(\text{e}_{6}+\sqrt{3}\text{e}_{4}\right).
\end{array}\label{eq:Tau(Octonions)}
\end{equation}
Such definition linearly extends to an order-three automorphism over
octonions $\mathbb{O}$ once we consider $\left\{ \text{e}_{0}=1,\text{e}_{2},...,\text{e}_{7}\right\} $
a basis for this algebra. It is interesting to note that in the octonions
there are two Argand planes, generated by $\left\{ \text{e}_{2},\text{e}_{5}\right\} $
and $\left\{ \text{e}_{4},\text{e}_{6}\right\} $ on which the automorphism
$\tau$ acts as the cubic root of unity $\frac{1}{2}\left(1+\sqrt{3}\text{i}\right)$.
Given the order-three automorphism $\tau$ over the Hurwitz algebra
of octonions $\left(\mathbb{O},\cdot,n\right)$ we can now obtain
a Petersson algebra \cite{Petersson 1969} defining a new product
\begin{equation}
x*y=\tau\left(\overline{x}\right)\cdot\tau^{2}\left(\overline{y}\right),
\end{equation}
for every $x,y\in\mathbb{O}$. The new algebra $\left(\mathbb{O},*,n\right)$
is again a composition algebra, called the Okubo algebra $\mathcal{O}$
(see \cite{Okubo1978b,Elduque Myung 90,Elduque 91,Elduque 93}). The
algebra does not have a unit nor a paraunit, but it has idempotent
elements. Again, it is worth noting that also the Okubo algebra is
a division algebra since if 
\begin{equation}
x*y=\tau\left(\overline{x}\right)\cdot\tau^{2}\left(\overline{y}\right)=0,
\end{equation}
then either $\tau\left(\overline{x}\right)$ or $\tau^{2}\left(\overline{y}\right)$
are zero and since $\tau$ is an automorphism, either $x$ or $y$
are zero. As in the paraoctonionic case if we use the split-octonions
as starting point we obtain $\left(\mathbb{O}_{s},*,n\right)$ which
is, again, a symmetric composition algebra but that has zero divisors
and it is thus called the split-Okubo algebra $\mathcal{O}_{s}$.

Unlike the case of the paraoctonions, it is worth presenting here
an independent and practical way of realising the Okubo algebra. In
fact, the following realisation was the one found by Okubo in \cite{Okubo 1978}.
Following \cite{Elduque Myung 90}, we define the real Okubo Algebra
$\mathcal{O}$ as the set of three by three Hermitian traceless matrices
over the complex numbers $\mathbb{C}$ with the following bilinear
product 
\begin{equation}
x*y=\mu\cdot xy+\overline{\mu}\cdot yx-\frac{1}{3}\text{Tr}\left(xy\right),\label{eq:product Ok}
\end{equation}
where $\mu=\nicefrac{1}{6}\left(3+\text{i}\sqrt{3}\right)$ and the
juxtaposition is the ordinary associative product between matrices.
It is worth noting that (\ref{eq:product Ok}) can be seen as a modification
of the Jordanian product which preserves the trace(lessness) of the
matrices. Indeed, setting $\mu=\nicefrac{1}{2}$ and disregarding
the last term, we retrieve the usual Jordan product over Hermitian
traceless matrices, i.e.
\begin{equation}
x\circ y=\frac{1}{2}xy+\frac{1}{2}yx.
\end{equation}
Nevertheless, Hermitian traceless matrices are not closed under such
product, thus requiring the additional term $-\nicefrac{1}{3}\text{Tr}\left(xy\right)$
for the closure of the algebra. Indeed, setting in (\ref{eq:product Ok})
$\text{Im}\mu=0$, one retrieves from the traceless part of the exceptional
Jordan algebra $\mathfrak{J}_{3}\left(\mathbb{C}\right)$, whose derivation
Lie algebra is $\mathfrak{su}\left(3\right)$.

Analyzing (\ref{eq:product Ok}), it becomes evident that the resulting
algebra is neither unital, associative, nor alternative. Nonetheless,
$\mathcal{O}$ is a\emph{ flexible }algebra, i.e. 
\begin{equation}
x*\left(y*x\right)=\left(x*y\right)*x,
\end{equation}
which will turn out to be an even more useful property than alternativity
in the definition of the projective plane. Even though the Okubo algebra
is not unital, it does have idempotents, i.e. $e*e=e$, such as 
\begin{equation}
e=\text{i}_{0}=\left(\begin{array}{ccc}
2 & 0 & 0\\
0 & -1 & 0\\
0 & 0 & -1
\end{array}\right),\label{eq:idemp}
\end{equation}
that together with
\begin{equation}
\begin{array}{ccc}
\text{i}_{1}=\sqrt{3}\left(\begin{array}{ccc}
0 & 1 & 0\\
1 & 0 & 0\\
0 & 0 & 0
\end{array}\right), &  & \text{i}_{2}=\sqrt{3}\left(\begin{array}{ccc}
0 & 0 & 1\\
0 & 0 & 0\\
1 & 0 & 0
\end{array}\right),\\
\text{i}_{3}=\sqrt{3}\left(\begin{array}{ccc}
0 & 0 & 0\\
0 & 0 & 1\\
0 & 1 & 0
\end{array}\right), &  & \text{i}_{4}=\sqrt{3}\left(\begin{array}{ccc}
1 & 0 & 0\\
0 & -1 & 0\\
0 & 0 & 0
\end{array}\right),\\
\text{i}_{5}=\sqrt{3}\left(\begin{array}{ccc}
0 & -i & 0\\
i & 0 & 0\\
0 & 0 & 0
\end{array}\right), &  & \text{i}_{6}=\sqrt{3}\left(\begin{array}{ccc}
0 & 0 & -i\\
0 & 0 & 0\\
i & 0 & 0
\end{array}\right),\\
\text{i}_{7}=\sqrt{3}\left(\begin{array}{ccc}
0 & 0 & 0\\
0 & 0 & -i\\
0 & i & 0
\end{array}\right),
\end{array}\label{eq:definizione i ottonioniche}
\end{equation}
form a basis for $\mathcal{O}$ that has real dimension $8$. It is
worth noting that the choice of the idempotent $e$ as in (\ref{eq:idemp})
does not yield to any loss of generality for the subsequent development
of our work since all idempotents are conjugate under the automorphism
group (cfr. \cite[Thm. 20]{ElduQueAut}). The choice of this special
basis is motivated on the fact that it will turn to be an orthonormal
basis with respect to the norm in (\ref{eq:Norm-Ok}) and that changing
the product over Okubo algebra according to (\ref{eq:Ok->Oct}) the
elements of the basis $\left\{ e=\text{i}_{0},\text{i}_{1},...,\text{i}_{7}\right\} $
will correspond to the octonionic one previously defined. 

In this realisation a direct way of defining a quadratic norm $n$
over Okubo algebra, is the following 
\begin{equation}
n\left(x\right)=\frac{1}{6}\text{Tr}\left(x^{2}\right),\label{eq:Norm-Ok}
\end{equation}
where $x^{2}$ is square of the element $x$ through the standard
matrix product for every $x\in\mathcal{O}$. It is straightforward
to see that the \emph{norm} has signature $(8,0)$, is associative
and composition over the Okubo algebra itself.

\subsection{\label{subsec:Okubo-algebra,-octonions}Okubo algebra, octonions
and paraoctonions}

Octonions, paraoctonions and Okubo algebras are mutually interconnected
in such a way that we can easily switch from one to the other by simply
changing the definition of the bilinear product over the vector space
of the algebra. Let us consider a new product over the Okubo algebra
$\mathcal{O}$ as
\begin{equation}
x\cdot y=\left(e*x\right)*\left(y*e\right),\label{eq:Ok->Oct}
\end{equation}
where $x,y\in\mathcal{O}$ and $e$ is an idempotent of $\mathcal{O}$.
Given that $e*e=e$ and $n\left(e\right)=1$, the element $e$ acts
as a left and right identity, i.e. 
\begin{align}
x\cdot e & =e*x*e=n\left(e\right)x=x,\\
e\cdot x & =e*x*e=n\left(e\right)x=x.
\end{align}
Moreover, since Okubo algebra is a composition algebra, the same norm
$n$ enjoys the following relation 
\begin{equation}
n\left(x\cdot y\right)=n\left(\left(e*x\right)*\left(y*e\right)\right)=n\left(x\right)n\left(y\right),
\end{equation}
which means that $\left(\mathcal{O},\cdot,n\right)$ is a unital composition
algebra of real dimension $8$. Since it is also a division algebra,
then it must be isomorphic to the octonions $\mathbb{O}$, as noted
by Okubo himself \cite{Okubo 1978,Okubo 78c}. On the other hand,
as already noticed, if we consider the order-three automorphism of
the octonions in (\ref{eq:Tau(Octonions)}), the Okubo algebra is
then realised as a Petersson algebra from the octonions setting
\begin{align}
x*y & =\tau\left(\overline{x}\right)\cdot\tau^{2}\left(\overline{y}\right).\label{eq:Oct->Ok}
\end{align}
Note that (\ref{eq:Tau(Octonions)}) is formulated assuming the knowledge
of the octonionic product. Reading the same maps as Okubic maps we
then have the notable relations (see \cite{ElduQueAut}), i.e.

\begin{align}
\overline{x} & =\left(\left(x*e\right)*e\right)*e,\\
\tau\left(x\right) & =\left(\left(\left(x*e\right)*e\right)*e\right)*e,
\end{align}
so that, in fact, the two maps are linked one another, i.e.

\begin{align}
\tau\left(x\right) & =\overline{x}*e\\
\overline{x} & =\tau\left(e*x\right).
\end{align}
While the order-three map $\tau$ and the octonionic conjugation are
intertwined, it is important to highlight their distinct impacts on
the algebras $\mathbb{O}$ and $\mathcal{O}$. While $\tau$ is an
automorphism for both Okubo algebra $\mathcal{O}$ and octonions $\mathbb{O}$,
$\overline{x}$ does not respect the algebraic structure of the Okubo
algebra $\mathcal{O}$, since it is not an automorphism nor an anti-automorphism
with respect to the Okubic product, while it is an anti-homomorphism
over octonions $\mathbb{O}$.

The scenario with paraoctonions, $p\mathbb{O}$ is more straightforward.
By definition, paraoctonions are obtainable from octonions $\mathbb{O}$,
through
\begin{equation}
x\bullet y=\overline{x}\cdot\overline{y},\label{eq:Oct->pOct}
\end{equation}
while, on the other hand, octonions $\mathbb{O}$ are obtainable from
paraoctonions $p\mathbb{O}$ through the aid of the para-unit $\boldsymbol{1}\in p\mathbb{O}$,
such that
\begin{align}
x\cdot y & =\left(\boldsymbol{1}\bullet x\right)\bullet\left(y\bullet\boldsymbol{1}\right)=\overline{x}\bullet\overline{y}.\label{eq:pOct->Oct}
\end{align}
The new algebra $\left(p\mathbb{O},\cdot,n\right)$ is again an eight-dimensional
composition algebra which is also unital and division and thus, for
Hurwitz theorem, isomorphic to the octonions $\mathbb{O}$. Moreover,
since $\tau\left(\overline{x}\right)=\overline{\tau\left(x\right)}$,
we also have that the Okubic algebra is obtainable from the para-Hurwitz
algebra with the introduction of a Petersson-like product, i.e. 
\begin{equation}
x*y=\tau\left(x\right)\bullet\tau^{2}\left(y\right).\label{eq:pOct->Ok}
\end{equation}
We thus have shown that all three algebras are obtainable one from
the other. 

A similar setup occurs for the split versions of the algebras. Indeed,
defining the product in (\ref{eq:Ok->Oct}) from $x,y\in\mathcal{O}_{s}$
and the product being the split-Okubic product, we again obtain a
unital composition algebra of dimension 8, but this time the algebra
is not a division algebra since the split-Okubo algebra $\mathcal{O}_{s}$
has zero divisors. Therefore, for the generalised Hurwitz theorem
the algebra $\left(\mathcal{O}_{s},\cdot,n\right)$ must be isomorphic
to that of the split-octonions $\mathbb{O}_{s}$. The same reasoning
can be applied to the product in (\ref{eq:pOct->Oct}) with $x,y\in p\mathbb{O}_{s}$
which again results in an unital composition algebra of dimension
8 which is not a division algebra and thus isomorphic to that of the
split-octonions $\mathbb{O}_{s}$. A final summary of all the relations
between the products of the 8-dimensional composition algebras is
in Table \ref{tab:Oku-Para-Octo}. 
\begin{table}
\begin{centering}
\begin{tabular}{|c|c|c|c|}
\hline 
Algebras & $\left(\mathcal{O},*,n\right)$ & $\left(p\mathbb{O},\bullet,n\right)$ & $\left(\mathbb{O},\cdot,n\right)$\tabularnewline
\hline 
\hline 
$x*y$ & $x*y$ & $\tau\left(x\right)\bullet\tau^{2}\left(y\right)$ & $\tau\left(\overline{x}\right)\cdot\tau^{2}\left(\overline{y}\right)$\tabularnewline
\hline 
$x\bullet y$ & $\tau^{2}\left(x\right)*\tau\left(y\right)$ & $x\bullet y$ & $\overline{x}\cdot\overline{y}$\tabularnewline
\hline 
$x\cdot y$ & $\left(e*x\right)*\left(y*e\right)$ & $\left(\boldsymbol{1}\bullet x\right)\bullet\left(y\bullet\boldsymbol{1}\right)$ & $x\cdot y$\tabularnewline
\hline 
\end{tabular}
\par\end{centering}
\medskip{}

\centering{}%
\begin{tabular}{|c|c|c|c|}
\hline 
split-Algebras & $\left(\mathcal{O}_{s},*,n\right)$ & $\left(p\mathbb{O}_{s},\bullet,n\right)$ & $\left(\mathbb{O}_{s},\cdot,n\right)$\tabularnewline
\hline 
\hline 
$x*y$ & $x*y$ & $\tau\left(x\right)\bullet\tau^{2}\left(y\right)$ & $\tau\left(\overline{x}\right)\cdot\tau^{2}\left(\overline{y}\right)$\tabularnewline
\hline 
$x\bullet y$ & $\tau^{2}\left(x\right)*\tau\left(y\right)$ & $x\bullet y$ & $\overline{x}\cdot\overline{y}$\tabularnewline
\hline 
$x\cdot y$ & $\left(e*x\right)*\left(y*e\right)$ & $\left(\boldsymbol{1}\bullet x\right)\bullet\left(y\bullet\boldsymbol{1}\right)$ & $x\cdot y$\tabularnewline
\hline 
\end{tabular}\caption{\label{tab:Oku-Para-Octo}\emph{On the left}: we show to obtain the
Okubic product $*,$ the paraoctonionic product $\bullet$ and the
octonionic product $\cdot$ from the Okubo algebra $\left(\mathcal{O},*,n\right)$,
paraoctonions $\left(p\mathbb{O},\bullet,n\right)$ and octonions
$\left(\mathbb{O},\cdot,n\right)$ respectively. \emph{On the right}:
we show the relations between the corresponding split-algebras.}
\end{table}

Finally, should be stressed that while transitioning from one algebra
to another is feasible, none of these algebras is isomorphic to another.
For example, while the octonions $\mathbb{O}$ are alternative and
unital, paraoctonions $p\mathbb{O}$ are neither alternative nor unital
but they do have a para-unit. In contrast, the Okubo algebra $\mathcal{O}$
is non-alternative, but rather it contains idempotent elements.

\subsection{Planes over symmetric composition algebras}

We now proceed as in Section \ref{sec:Veronese-formulation-of} to
define in a general setting all projective and hyperbolic planes over
the paraunital and non-unital algebras $p\mathbb{O},p\mathbb{O}_{s},\mathcal{O},\mathcal{O}_{s}$.
Let $\widetilde{V}=\mathbb{K}^{3}\times\mathbb{R}^{3}$ where $\mathbb{K}\in\left\{ p\mathbb{O},p\mathbb{O}_{s},\mathcal{O},\mathcal{O}_{s}\right\} $.
For a triple $\left(\gamma_{1},\gamma_{2},\gamma_{3}\right)$ with
$\gamma_{1},\gamma_{2},\gamma_{3}\in\left\{ \pm1\right\} $, we define
a Veronese-like vector as a non-zero vector $\left(x_{\nu};\lambda_{\nu}\right)_{\nu}\in\widetilde{V}$
such that 
\begin{align}
\lambda_{1}x_{1} & =\gamma_{3}\gamma_{2}^{-1}x_{2}\circ x_{3},\,\,\lambda_{2}x_{2}=\gamma_{1}\gamma_{3}^{-1}x_{3}\circ x_{1},\,\,\,\lambda_{3}x_{3}=\gamma_{2}\gamma_{1}^{-1}x_{1}\circ x_{2},\label{eq:VerGen1-1-1}\\
n\left(x_{1}\right) & =\gamma_{3}\gamma_{2}^{-1}\lambda_{2}\lambda_{3},\,\,n\left(x_{2}\right)=\gamma_{1}\gamma_{3}^{-1}\lambda_{3}\lambda_{1},\,\,\,n\left(x_{3}\right)=\gamma_{2}\gamma_{1}^{-1}\lambda_{1}\lambda_{2},\label{eq:VerGen2-1-1}
\end{align}
where $x_{\nu}\in\mathbb{K}$ and $\lambda_{\nu}\in\mathbb{R}$ for
$\nu=1,2,3$ and $\circ$ represents the appropriate Okubic, split-Okubic,
paraoctonionic or split-paraoctonionic product. We call such subset
$\widetilde{H}$. \emph{Points} $\mathscr{P}$ of the corresponding
incidence plane are defined as the 1-dimensional subspaces $\mathbb{R}w$,
where $w$ is a Veronese-like vector, i.e.,
\begin{equation}
\mathscr{P}=\left\{ \mathbb{R}w:w\in\widetilde{H}\smallsetminus\left\{ 0\right\} \right\} .\label{eq:Octonionic point-1-1}
\end{equation}
 As for the octonionic line we proceed defining the bilinear form
$\beta$ as
\begin{equation}
\beta\left(x,y\right)=\stackrel[\nu=1]{3}{\sum}\left(\gamma_{\nu+2}^{-1}\gamma_{\nu+1}\left\langle x_{\nu},y_{\nu}\right\rangle +\lambda_{\nu}\mu_{\nu}\right),
\end{equation}
where $x=\left(x_{\nu};\lambda_{\nu}\right)_{\nu}$ and $y=\left(y_{\nu};\mu_{\nu}\right)_{\nu}$
and indices are taken modulo $3$. It is worth noting that while the
product changes between the 8-dimensional composition algebras, their
norm $n$ (and thus its polarisation) remains the same. As in previous
constructions, the set of the lines $\mathscr{L}$ are given as the
orthogonal space to a Veronese vector $x$, i.e.

\begin{equation}
\ell\coloneqq x^{\perp}=\left\{ y\in\mathbb{K}^{3}\times\mathbb{R}^{3}:\beta\left(y,x\right)=0\right\} ,
\end{equation}
which explicitly happens if and only if

\begin{equation}
\gamma_{3}^{-1}\gamma_{2}\left\langle x_{1},y_{1}\right\rangle +\gamma_{1}^{-1}\gamma_{3}\left\langle x_{2},y_{2}\right\rangle +\gamma_{2}^{-1}\gamma_{1}\left\langle x_{3},y_{3}\right\rangle +\lambda_{1}\mu_{1}+\lambda_{2}\mu_{2}+\lambda_{3}\mu_{3}=0.\label{eq:Explicit-line-1-1}
\end{equation}
Finally, we remark that the incidence relations remain unchanged,
they are always determined by inclusion. 

\subsection{Isomorphisms with the octonionic planes}

Given the previous definitions, for every $\mathbb{K}\in\left\{ p\mathbb{O},p\mathbb{O}_{s},\mathcal{O},\mathcal{O}_{s}\right\} $
we then have a projective plane $\mathbb{K}P^{2}$ for $\left(\gamma_{1},\gamma_{2},\gamma_{3}\right)=\left(1,1,1\right)$
and a hyperbolic plane $\mathbb{K}H^{2}$ setting $\left(\gamma_{1},\gamma_{2},\gamma_{3}\right)=\left(1,1,-1\right)$.
In \cite{CMZ-2} we have proved that this definition of the paraoctonionic
projective plane $p\mathbb{O}P^{2}$ and the Okubic projective plane
$\mathcal{O}P^{2}$ are isomorphic to the octonionic projective plane
$\mathbb{O}P^{2}$. Now, it is a straightforward calculation to prove
the isomorphism between $\mathcal{O}_{s}P^{2}$ and $\mathbb{O}_{s}P^{2}$,
realised by the same map involving the division counterparts,

\begin{equation}
\Phi:\begin{cases}
\left(x,y,x*y;n\left(y\right),n\left(x\right),1\right)\longrightarrow\left(\tau^{2}\left(\overline{x}\right),y,y\cdot\overline{\tau^{2}\left(\overline{x}\right)};n\left(y\right),n\left(x\right),1\right),\\
\left(0,0,x;n\left(x\right),1,0\right)\longrightarrow\left(0,0,\tau^{2}\left(\overline{x}\right);n\left(x\right),1,0\right),\\
\left(0,0,0;1,0,0\right)\longrightarrow\left(0,0,0;1,0,0\right),
\end{cases}\label{eq:split-OkuboIsomorph}
\end{equation}
where the left-hand side vectors are Veronese under conditions (\ref{eq:VerGen1-1-1})
and (\ref{eq:VerGen2-1-1}) with $\mathbb{K}=\mathcal{O}_{s}$ and
$\left(\gamma_{1},\gamma_{2},\gamma_{3}\right)=\left(1,1,1\right)$,
while the right-hand side vectors are Veronese under conditions (\ref{eq:VerGen1-1})
and (\ref{eq:VerGen2-1}) with $\mathbb{K}=\mathbb{O}_{s}$ and $\left(\gamma_{1},\gamma_{2},\gamma_{3}\right)=\left(1,1,1\right)$.

A similar map establishes the isomorphism between $p\mathbb{O}_{s}P^{2}$
and $\mathbb{O}_{s}P^{2}$. Indeed, this isomorphism is realised by 

\begin{equation}
p\Phi:\begin{cases}
\left(x,y,x\bullet y;n\left(y\right),n\left(x\right),1\right)\longrightarrow\left(\tau^{2}\left(x\right),y,y\cdot\overline{\tau^{2}\left(x\right)};n\left(y\right),n\left(x\right),1\right),\\
\left(0,0,x;n\left(x\right),1,0\right)\longrightarrow\left(0,0,\tau^{2}\left(x\right);n\left(x\right),1,0\right),\\
\left(0,0,0;1,0,0\right)\longrightarrow\left(0,0,0;1,0,0\right),
\end{cases}\label{eq:paraOctonions-isomorp}
\end{equation}
where the left-hand side vectors are Veronese under conditions (\ref{eq:VerGen1-1-1})
and (\ref{eq:VerGen2-1-1}) with $\mathbb{K}=p\mathbb{O}_{s}$ and
$\left(\gamma_{1},\gamma_{2},\gamma_{3}\right)=\left(1,1,1\right)$,
the right-hand side vectors are Veronese under conditions (\ref{eq:VerGen1-1})
and (\ref{eq:VerGen2-1}) with $\mathbb{K}=\mathbb{O}_{s}$ and $\left(\gamma_{1},\gamma_{2},\gamma_{3}\right)=\left(1,1,1\right)$.
Finally, a careful inspection of the (\ref{eq:split-OkuboIsomorph})
and (\ref{eq:paraOctonions-isomorp}) shows that the two maps $\Phi$
and $p\Phi$ are also isomorphisms between the Okubic hyperbolic plane
$\mathcal{O}H^{2}$, the paraoctonionic hyperbolic plane $p\mathbb{O}H^{2}$
and the octonionic plane $\mathbb{O}H^{2}$ respectively, when one
sets $\left(\gamma_{1},\gamma_{2},\gamma_{3}\right)=\left(1,1,-1\right)$
in both the generalised Veronese conditions (\ref{eq:VerGen1-1-1})-(\ref{eq:VerGen2-1-1})
and (\ref{eq:VerGen1-1})-(\ref{eq:VerGen2-1}). 

In summary, the following isomorphisms among planes over 8-dimensional
symmetric composition algebras hold:
\begin{align}
\mathcal{O}P^{2} & \cong p\mathbb{O}P^{2}\cong\mathbb{O}P^{2},\\
\mathcal{O}_{s}P^{2} & \cong p\mathbb{O}_{s}P^{2}\cong\mathbb{O}_{s}P^{2},\\
\mathcal{O}H^{2} & \cong p\mathbb{O}H^{2}\cong\mathbb{O}H^{2}.
\end{align}

\section{\label{sec:Complex-Cayley-plane}Complex Cayley plane }

To conclude this paper, it is worth pointing out that also the complex
Lie groups $\text{E}_{6}^{\mathbb{C}}$,$\text{F}_{4}^{\mathbb{C}}$
and $\text{G}_{2}^{\mathbb{C}}$ respectively, can be realised geometrically
as the collineation group, the isometry group and the stabiliser of
a quadrangle of a projective plane realised within the same Veronese
framework. Indeed, all constructions we have previously done can be
easily generalised \emph{for any field}. More specifically, when we
considered the complex octonions $\mathbb{O}_{\mathbb{C}}$ defined
in Table \ref{tab:Split Octonions-1-1}, or, equivalently, as $\mathbb{O}_{\mathbb{C}}\cong\mathbb{\mathbb{C}\otimes\mathbb{O}}\cong\mathbb{\mathbb{C}\otimes\mathbb{O}}_{s}$,
and one can set a complex norm $n_{\mathbb{C}}$ from complex octonions
$\mathbb{O}_{\mathbb{C}}$ to the complex field $\mathbb{C}$ as 
\begin{equation}
n_{\mathbb{C}}\left(b\right)=z_{0}^{2}+...+z_{7}^{2}\in\mathbb{C},
\end{equation}
where $b=\stackrel[k=0]{7}{\sum}z_{k}\text{i}_{k}$ with $z_{k}\in\mathbb{C}$.
It is straightforward to see that $\left(\mathbb{O}_{\mathbb{C}},\cdot,n_{\mathbb{C}}\right)$
is a composition algebra, i.e. $n_{\mathbb{C}}\left(b_{1}b_{2}\right)=n_{\mathbb{C}}\left(b_{1}\right)n_{\mathbb{C}}\left(b_{2}\right)$.
Moreover, the involution $b\longrightarrow\overline{b}$ with
\begin{equation}
\overline{b}=z_{0}\text{i}_{0}-\stackrel[k=1]{7}{\sum}z_{k}\text{i}_{k},
\end{equation}
is such that $n_{\mathbb{C}}\left(b\right)=b\overline{b}$. One can
then define the complex Cayley plane from the Veronese conditions
with the following construction. Let $V_{\mathbb{C}}\cong\mathbb{\mathbb{O}_{\mathbb{C}}}^{3}\times\mathbb{C}^{3}$
be a complex vector space, with elements $\omega$ of the form 
\[
\left(b_{\nu};\lambda_{\nu}\right)_{\nu}=\left(b_{1},b_{2},b_{3};\lambda_{1},\lambda_{2},\lambda_{3}\right),
\]
where $b_{\nu}\in\mathbb{\mathbb{C}\otimes\mathbb{O}}$, $\lambda_{\nu}\in\mathbb{C}$
for $\nu=1,2,3$. A Veronese vector $\omega\in V_{\mathbb{C}}$ is
now given by the following conditions: 

\begin{align}
\lambda_{1}\overline{b}_{1} & =b_{2}b_{3},\,\,\lambda_{2}\overline{b}_{2}=b_{3}b_{1},\,\,\lambda_{3}\overline{b}_{3}=b_{1}b_{2},\label{eq:Veronese conditions1}\\
n_{\mathbb{C}}\left(b_{1}\right) & =\lambda_{2}\lambda_{3},\,n_{\mathbb{C}}\left(b_{2}\right)=\lambda_{3}\lambda_{1},n_{\mathbb{C}}\left(b_{1}\right)=\lambda_{1}\lambda_{2}.\nonumber 
\end{align}
Let the set $H\subset V_{\mathbb{C}}$ be the set of Veronese vectors
inside $V_{\mathbb{C}}$. Since $\mathbb{C}$ is commutative, and
$\lambda b=b\lambda$, and $n_{\mathbb{C}}\left(\lambda b\right)=\lambda^{2}n_{\mathbb{C}}\left(b\right)$
when $\lambda\in\mathbb{C}$, then if $\omega$ is a Veronese vector,
all complex multiples $\mathbb{C\omega}$ are again Veronese vectors,
i.e. if $\omega\in H$ then $\mu\omega\in H$ when $\mu\in\mathbb{C}$.
The \emph{complex Cayley plane} $\mathbb{O}P_{\mathbb{C}}^{2}$ has
the set of 1-dimensional complex subspaces $\mathbb{C}\omega$ as
\emph{points} of the plane, i.e.
\begin{equation}
\mathbb{O}P_{\mathbb{C}}^{2}=\left\{ \mathbb{C}\omega:\omega\in H\smallsetminus\left\{ 0\right\} \right\} ,
\end{equation}
and as \emph{line} $\ell$ has the orthogonal subspace 
\begin{equation}
\ell\coloneqq\omega^{\perp}=\left\{ \upsilon\in V_{\mathbb{C}}:\beta\left(\upsilon,\omega\right)=0\right\} ,
\end{equation}
where the complex bilinear form $\beta$ is defined as

\begin{equation}
\beta\left(\upsilon,\omega\right)=\stackrel[\nu=1]{3}{\sum}\left(\left\langle b_{\nu}^{1},b_{\nu}^{2}\right\rangle _{\mathbb{O}_{\mathbb{C}}}+\lambda_{\nu}^{1}\lambda_{\nu}^{2}\right),
\end{equation}
with $\upsilon,\omega\in V_{\mathbb{C}}$, of coordinates $\left(b_{\nu}^{1};\lambda_{\nu}^{1}\right)_{\nu}$,$\left(b_{\nu}^{2};\lambda_{\nu}^{2}\right)_{\nu}$
respectively. 

Again, the map (\ref{eq:mappingVer-Jordan}) from $V_{\mathbb{C}}$
to the complexification of the exceptional Jordan algebra $\mathfrak{J}_{3}^{\mathbb{C}}\left(\mathbb{O}\right)\cong\mathbb{C}\otimes\mathfrak{J}_{3}\left(\mathbb{O}\right)$
(namely, to the complex Albert algebra) establishes the one-to-one
correspondence between Veronese vectors and rank-1 idempotent elements
of $\mathfrak{J}_{3}^{\mathbb{C}}\left(\mathbb{O}\right)$. It is
thus straightforward to state the alternative definition of the complex
Cayley-Moufang plane as

\begin{equation}
\mathbb{O}_{\mathbb{C}}P^{2}\cong\left\{ X\in\mathfrak{J}_{3}^{\mathbb{C}}\left(\mathbb{O}\right):X^{\#}=0,\text{tr}\left(X\right)=1\right\} ,
\end{equation}
and the consequent characterization of the collineation groups as
shown in Table \ref{tab:Relations-between-incidence}. Non-alternative
constructions for $\mathcal{O}_{\mathbb{C}}P^{2}$ and $p\mathbb{O}_{\mathbb{C}}P^{2}$
as in Section \ref{sec:Realisations-through-symmetric} along with
their isomorphisms with $\mathbb{O}_{\mathbb{C}}P^{2}$ the are straightforward
and are left to the reader.
\begin{center}
\begin{figure}
\begin{centering}
\includegraphics[scale=0.7]{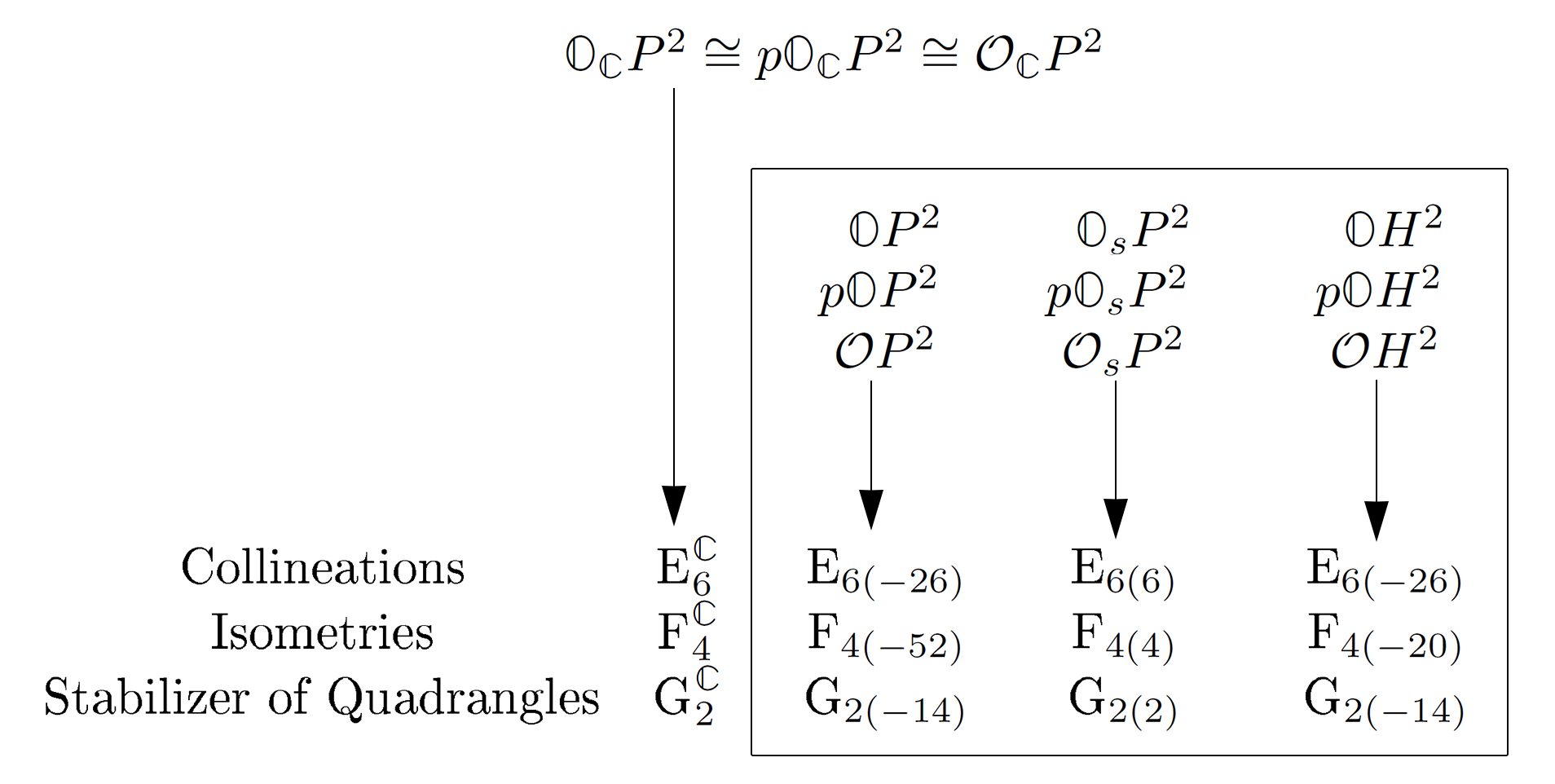}\caption{\label{fig:Collineation-groups}Collineation groups}
\par\end{centering}
\end{figure}
\par\end{center}

\section{Conclusions e further developments}

In this article, we presented a comprehensive framework for constructing
both octonionic and split-octonionic projective and hyperbolic planes
(in their Riemannian and pseudo-Riemannian forms), exploiting the
formalism of Veronese vectors. This unified approach provides a concrete
and efficient method for representing these incidence planes. Furthermore,
this approach has the advantage to manifestly exhibit the relation
between the aforementioned incidence planes and the exceptional Lie
groups of type $\text{E}_{6}$, $\text{F}_{6}$, and $\text{G}_{2}$,
which are then respectively realised as collineation groups, isometry
groups, and stabilizers of quadrangles of the aforementioned planes.
Also, we introduced novel realizations of the hyperbolic octonionic
plane and split-octonionic projective plane using the non-unital,
non-alternative symmetric composition algebras such as the (split-)Okubo
algebra and the (split-)paraoctonionic algebra. Not only these constructions
have a mathematical relevance on their own, but they can also pertain
to interesting physical applications. 

Finally, we used the same Veronese framework for the definition of
the complex Cayley plane which is linked to the complexification of
the exceptional Jordan algebra (\emph{a.k.a.} Albert algebra), $\mathfrak{J}_{3}^{\mathbb{C}}\left(\mathbb{O}\right)$
. This naturally yields to the complex forms of the Lie groups $\text{E}_{6},\text{F}_{4}$
and $\text{G}_{2}$. This result completes the geometrical and algebraic
landscape depicted in Figure \ref{fig:Collineation-groups}. 

Concluding, let us point out that not only our framework provides
a constructive and sound mathematical approach to the various forms
of incidence planes over eight-dimensional division algebras, but
it also paves the way to the understanding of the geometries involving
exceptional Lie groups. It is our hope that our work may give rise
to new geometric perspectives in areas of Mathematics and Physics
where exceptional Lie groups are of paramount importance.

\section{Acknowledgments}

The work of AM is supported by a \textquotedblleft Maria Zambrano\textquotedblright{}
distinguished researcher fellowship, financed by the European Union
within the NextGenerationEU program.

$*$\noun{ Departamento de Matemática, }\\
\noun{Universidade do Algarve, }\\
\noun{Campus de Gambelas, }\\
\noun{8005-139 Faro, Portugal} 
\begin{verbatim}
a55499@ualg.pt

\end{verbatim}
$\dagger$\noun{ Instituto de Física Teorica, Dep.to de Física,}\\
\noun{Universidad de Murcia, }\\
\noun{Campus de Espinardo, }\\
\noun{E-30100, Spain}
\begin{verbatim}
alessio.marrani@um.es 

\end{verbatim}
$\ddagger$\noun{ Dipartimento di Scienze Matematiche, Informatiche
e Fisiche, }\\
\noun{Università di Udine, }\\
\noun{Udine, 33100, Italy} 
\begin{verbatim}
francesco.zucconi@uniud.it
\end{verbatim}

\begin{verbatim}

\end{verbatim}

\end{document}